\documentclass{tran-l2}
\usepackage{graphicx}
\usepackage{amssymb}

\newtheorem{proposition}{Proposition}[section]
\newtheorem{theorem}{Theorem}[section]
\newtheorem{lemma}[theorem]{Lemma}
\newtheorem{cor}[theorem]{Corollary}

\theoremstyle{definition}

\theoremstyle{remark}

\numberwithin{equation}{section}

\newcommand{\A }{\mathcal{A}}

\newcommand{\B }{\mathcal{B}}

\newcommand{\ch }{\mathrm{ch}}
\newcommand{\cD}{\mathcal{D}}
\newcommand{\cA }{\mathcal{A}}
\newcommand{\cI}{\mathcal{I}}
\newcommand{\cL}{\mathcal{L}}

\newcommand{\cO}{\mathcal{O}}

\newcommand{\cV}{\mathcal{V}}
\newcommand{\cX}{\mathcal{X}}
\newcommand{\cbX}{\overline{\mathcal{X}}}

\newcommand{\C }{\mathbb{C}}
\newcommand{\dif }{\mathrm{d}}
\newcommand{\e }{\mathrm{e}}

\newcommand{\Gam }{\varGamma }
\newcommand{\gamgr }{\Gamma^1_q(\mathrm{Gr}(r,N))}
\newcommand{\gamgrez }{\Gamma^1_q(\mathrm{Gr}(1,2))}
\newcommand{\gfrak}{\mathfrak{g}}
\newcommand{\Gr}{\mbox{Gr}(r,N)}

\newcommand{\Hom}{\mathrm{Hom}}
\newcommand{\id}{\mbox{Id}}
\newcommand{\im}{\mbox{Im}}
\newcommand{\kopr }{\varDelta }
\newcommand{\kow }{\varDelta }
\newcommand{\K }{K}

\newcommand{\Lin }{\mathrm{Lin}}
\newcommand{\mat}{\mathrm{Mat}}
\newcommand{\N}{\mathbb{N}}
\newcommand{\ogr}{\mathcal{O}(\mathrm{Gr}(r,N))}
\newcommand{\ot }{\otimes }
\newcommand{\opartial }{\overline{\partial} }

\newcommand{\oplr}{\mathcal{O}_q({\mathcal V}_r)}

\newcommand{\oqg}{{\mathcal O}_q(G)}
\newcommand{\oqsl}{{\mathcal O}_q(\mathrm{SL}(N))}
\newcommand{\oqsu}{{\mathcal O}_q(\mathrm{SU}(N))}
\newcommand{\oqgr }{\mathcal{O}_q(\mathrm{Gr}(r,N))}
\newcommand{\oqgrez }{\mathcal{O}_q(\mathrm{Gr}(1,2))}
\newcommand{\oqgrminus }{\mathcal{O}_q(\mathrm{Gr}(r{-}1,N{-}1))}

\newcommand{\ovr}{{\mathcal O}({\mathcal V}_r)}

\newcommand{\oqvr}{{\mathcal O}_q({\mathcal V}_r)}
\newcommand{\pair }[2]{\langle #1,#2\rangle }

\newcommand{\plr}{{\mathcal V}_r}

\newcommand{\qd }{\hat{q}}

\newcommand{\rf }{{\bf r}}
\newcommand{\rh }{\hat{R}}
\newcommand{\ra }{\acute{R}}
\newcommand{\rg }{\grave{R}}
\newcommand{\rc }{\check{R}}
\newcommand{\R }{\mathbb{R}}

\newcommand{\slfrak}{\mathfrak{sl}}
\newcommand{\sr }{ s }
\newcommand{\SL }{\mathrm{SL}(N)}
\newcommand{\SLN }{\mathcal{O}_q(\mathrm{SL}(N))}

\newcommand{\SU }{\mathrm{SU}(N)}
\newcommand{\SUN }{\mathcal{O}_q(\mathrm{SU}(N))}
\newcommand{\Sym }{\mathrm{Sym}}
\newcommand{\tr }{\mathrm{tr}}

\newcommand{\ub }{{\bf u}}

\newcommand{\ug }{U_q(\mathfrak{g})}
\newcommand{\UcslN }{U(\mathfrak{sl}_N)}
\newcommand{\UslN }{U_q(\mathfrak{sl}_N)}

\newcommand{\vep }{\varepsilon }

\begin{document}

\title[Differential Calculus on Quantum Grassmannians I]
 {Differential Calculus on Quantum Complex Grassmann Manifolds I:
 Construction}

\author{Stefan Kolb}
\address{Mathematisches Institut, Universit\"at Leipzig, Augustusplatz 10,
         04109 Leipzig, Germany}
\email{kolb@itp.uni-leipzig.de}
\thanks{This work was supported by the Deutsche Forschungsgemeinschaft
  within the scope of the postgraduate scholarship programme
  ``Graduiertenkolleg Quantenfeldtheorie'' at the
  University of Leipzig}

%    General info
%\subjclass[2000]{Primary 58B32, 81R50; Secondary 14M15}


\begin{abstract}
  Covariant first order differential calculus over quantum complex Grassmann
manifolds is considered. It is shown by a Pusz--Woronowicz type argument that
under restriction to calculi close to classical K\"ahler differentials there 
exist exactly two such calculi for the homogeneous coordinate ring.
Complexification and localization procedures are used to induce covariant
first order differential calculi over quantum Grassmann manifolds.
It is shown that these differential calculi behave in many respects as their
classical counterparts. As an example the $q$-deformed Chern character of the
tautological bundle is constructed.
\end{abstract}

\maketitle

Covariant first order differential calculus is a concept first introduced by
S.~L.~Wo\-ronowicz \cite{a-Woro3}, \cite{a-Woro2}
to generalize the notion of differential form from commutative algebra
to quantum groups and quantum spaces. The task to find well behaved
analogues of differential forms for the noncommutative deformed coordinate
algebras appearing in the framework of quantum groups is still of
considerable interest \cite{a-herm98}, \cite{a-Schmue1}, \cite{a-SinVa}.
In \cite{a-PuWo89} W.~Pusz and 
S.L.~Woronowicz proved that for the quantum vector space of dimension $\ge 3$
there exist exactly two differential calculi freely generated by the
differentials of the 
generators. In \cite{a-Po92} P.~Podle\'{s} classified differential structures
on the quantum 2-sphere $S^2_{qc}$ which have certain properties similar to
classical differential forms. It turned out, that only in the so called
quantum subgroup case $c=0$ such a differential calculus exists and is then
uniquely determined. 

Podle\'{s} quantum 2-sphere is an example of a $q$-deformed Grassmann
manifold. The undeformed Grassmann manifold $\Gr$ of $r$-dimensional
subspaces in $\C^N$ is a projective algebraic variety. It is well
known \cite{b-CP94} that its homogeneous coordinate ring
${\mathcal O}(\mathcal{V}_r)$ can be $q$-deformed to a quantum space 
${\mathcal O}_q(\mathcal{V}_r).$ 
On the other hand M.~Noumi, M.~S.~Dijkhuizen and T.~Sugitani introduced a 
large class of quantum Grassmannians \cite{a-NDS97}.
Here only the quantum subgroup case is considered and will be denoted
by ${\mathcal O}_q(\Gr)$. The aim of this paper is to construct a
canonical covariant first order differential calculus over $\oqgr$ and
investigate its properties. 

The $q$-deformed coordinate rings $\oqvr$ and $\oqgr$ are closely related. 
More explicitly we consider the complexification of $\oqvr$, an algebra
obtained by adding complex conjugate elements. This complexification has
a distinguished invariant element, and $\oqgr$ can be obtained as a
subalgebra of the localization with respect to this element.
It is the main observation of our work that this construction allows an
analogue on the level of differential calculus. To achieve this the notions
of complexification and localization of first order differential calculus are
introduced. In generalization of the above mentioned classification result
of \cite{a-PuWo89} differential calculi over $\oqvr$ are classified. 
For this purpose the $U({\mathfrak sl}_{N})$-module structure
of classical K\"ahler differentials over ${\mathcal O}(\mathcal{V}_r)$
is analyzed. Under suitable additional assumptions all covariant
first order differential calculi with the same 
$U_q({\mathfrak sl}_{N})$-module structure are classified and it is seen 
that exactly two such calculi exist. 
Complexification and localization then lead to a calculus $\gamgr$ over 
${\mathcal O}_q(\Gr)$ which in the case of $S^2_{q0}$ is seen to be
isomorphic to the calculus constructed by Podle\'{s}. Thus a relationship 
between the calculus of Pusz and Woronowicz on the quantum plane and 
Podle\'{s}' calculus on $S^2_{q0}$ is established and generalized to 
arbitrary Grassmann manifolds.

The construction of $\gamgr$ given here allows to perform explicit
computations. As an example it is shown that as a left module
$\gamgr$ is generated by the differentials of the generators of $\oqgr$.
The dimension of $\gamgr$ in the sense of \cite{a-herm01} is estimated.

Considering $\oqgr$ in terms of the appropriate set of generators and
relations it is straightforward to construct a $q$-deformed analogue of the
module of sections of the tautological bundle over $\Gr$. This allows to
introduce the $q$-deformed Chern character of the tautological bundle.
It is a sum of closed differential forms in the universal
higher order differential calculus of $\gamgr$ which are seen to be central.
Thus it is to be expected that the cohomology ring is independent of the
deformation parameter. 

The ordering of the paper is as follows. Section \ref{Notconv} serves to
fix notations. In Section \ref{grass} quantum Grassmann manifolds are
recalled. The notion of complexification and localization are applied to the
homogeneous
coordinate ring $\oqvr$ establishing the relation between $\oqvr$ and
$\oqgr$. This part owes much to the detailed analysis of \cite{phd-stok}.
Explicit relations for different sets of generators of $\oqgr$ are
established. In a slight digression these relations are used to calculate the
kernel of the quantum analogue of the canonical map dual to the inclusion
$\textrm{Gr}(r{-}1,N{-}1)\hookrightarrow \textrm{Gr}(r,N)$.
The explicit relations between
different sets of generators involve a considerable number of $R$-matrices.
To simplify notations the graphical calculus of \cite{b-Turaev} is used.
In Appendix \ref{graphcal} this calculus is recalled and useful
simplifications of morphisms often met in our framework are provided.
The explicit construction of the canonical covariant first order
differential calculus is performed in Section \ref{construction}. First the
notations of localization and complexification for differential calculus on
$q$-spaces are introduced. Then under suitable additional assumptions all
differential calculi over the homogeneous coordinate ring $\oqvr$ are
classified. They are factored by their torsion submodules and then turn out
to be given by relations of $R$-matrix type. Localization and
complexification is applied to induce the desired calculus over
$\oqgr$.

The last section is devoted to the construction of the $q$-deformed Chern
character for the tautological bundle over $\Gr$.

\section{Notations and conventions}\label{Notconv}
Unless stated otherwise all notations and conventions coincide with those
used in \cite{b-KS}. Throughout this work $q^k\neq 1$ for all $k\in\N$ is
assumed.

The $q$-deformed universal enveloping algebra $\UslN$ is the complex
algebra generated by elements $E_i,F_i,K_i$ and $K_i^{-1}$, $i=1,\dots,N{-}1$
and relations
\begin{align*}
K_iK_i^{-1}&=K_i^{-1}K_i=1&& \\
K_iK_j&=K_jK_i&&\\
K_iE_j&=q^{a_{ij}}E_jK_i&&\\
K_iF_j&=q^{-a_{ij}}F_jK_i&&\\
E_iF_j-F_jE_i&=\delta_{ij}\frac{K_i-K_i^{-1}}{q-q^{-1}}\\
\sum_{k=0}^{1-a_{ij}}(-1)^k\left(\begin{array}{c}1-a_{ij}\\k \end{array}
        \right)_q& E_i^{1-a_{ij}-k}E_jE_i^k=0,&i&\neq j\\
\sum_{k=0}^{1-a_{ij}}(-1)^k\left(\begin{array}{c}1-a_{ij}\\k \end{array}
        \right)_q&F_i^{1-a_{ij}-k}F_jF_i^k=0,&i&\neq j.
\end{align*}  
Here
\begin{align*}
  a_{ij}=\begin{cases} 2&\textrm{if $i=j$,}\\
                      -1&\textrm{if $|i{-}j|=1,$}\\
                       0&\textrm{else}\end{cases}
\end{align*}                     
denotes the Cartan matrix of $\mathfrak{sl}_N$ and the
$q$-deformed binomial coefficient is defined by
\begin{align*}
{n\choose k}_q=\frac{[n][n{-}1]\dots[n{-}k{+}1]}{[1][2]\dots[k]}
\end{align*}
where $[x]=\frac{q^x-q^{-x}}{q-q^{-1}}$. 
The algebra $\UslN$ obtains a Hopf algebra structure by
\begin{eqnarray*}
\kopr K_i&=& K_i\otimes  K_i\\   
\kopr E_i&=& E_i\otimes K_i+1\otimes E_i\\
\kopr F_i&=& F_i\otimes 1 + K_i^{-1}\otimes F_i\\
\epsilon(K_i)&=&1,\quad \epsilon(E_i)=\epsilon(F_i)=0\\
S(K_i)&=&K_i^{-1},\quad S(E_i)=-E_iK_i^{-1},\quad S(F_i)=-K_iF_i.
\end{eqnarray*}
Among the irreducible finite dimensional representations of $\UslN$ are the
so called type 1 representations $V(\lambda)$ which are uniquely determined
by a highest weight $\lambda$. More explicitly if
$\omega_i,\, i=1,\dots,N{-}1$ denote the fundamental weights and if
$\lambda=\sum_i\lambda_i\omega_i\in\sum_i\N_0\omega_i$ is 
a dominant integral weight there exists $v\in V(\lambda)$  such that
\begin{align}\label{hwproperty}
  F_i v=0,\qquad  K_i^{-1}v=q^{\lambda_i}v.
\end{align}  

Dually the $q$-deformed algebra of regular functions $\oqsl$ is
the subalgebra of the dual Hopf algebra of $\UslN$ generated
by the matrix coefficients of the irreducible type 1 representations.
It is generated by the matrix coefficients $u^i_j$, $i,j=1,\dots,N$
of the vector representation.
By construction the type 1 representations of $\UslN$ coincide with
the right $\oqsl$-comodules. If $q\in\R$ the Hopf algebra
$\oqsl$ can be endowed with a $*$-structure by $(u^i_j)^*=S(u^j_i)$
and is then denoted by $\oqsu$.
The category of type 1 representations of $\UslN$ is a braided
monoidal category. There exists a universal $r$-form
\begin{align*}
\rf:\oqsl\otimes\oqsl\rightarrow \C
\end{align*}
such that the braiding $\rh_{VW}:V\otimes W\rightarrow W\otimes V$
can be written as
\begin{align*}
  v\otimes w\mapsto w_{(0)}\otimes v_{(0)}\rf(v_{(1)},w_{(1)})
\end{align*}  
where Sweedler notation is used.
The convolution inverse of $\rf$ is denoted by $\overline{\rf}$.
For any $a,b\in\oqsl$ the relations
\begin{align}
\rf(a_{(1)},b_{(1)})a_{(2)}b_{(2)}&=b_{(1)}a_{(1)}\rf(a_{(2)},b_{(2)})
\label{rform}\\
\overline{\rf}(a_{(1)},b_{(1)})b_{(2)}a_{(2)}&=a_{(1)}b_{(1)}\overline{\rf}
(a_{(2)},b_{(2)})\label{rform-}
\end{align}  
hold. In terms of the solution
\begin{align}\label{R}
\rh^{kl}_{ij}=\begin{cases}q&\textrm{if $i=j=k=l$,}\\
                           1&\textrm{if $i=l\neq j=k$,}\\
                           q-q^{-1}&\textrm{if $k=i<j=l$,}\\
                           0&\textrm{else}
              \end{cases}
\end{align}
of the braid relation the universal $r$-form is determined by
\begin{align}
  \rf(u^l_i,u^k_j)=p\rh_{ij}^{kl}\label{rpR}
\end{align}  
where $p\in\C$ satisfies $p^N=q^{-1}$. This universal $r$-form is real, i.e.
\begin{equation*}
  \rf(a^*,b^*)=\overline{\rf(b,a)}.
\end{equation*}  
Introduce variants of the $\rh$-matrix by
\begin{align*}
  \rc^{kl}_{ij}&=\rh^{ji}_{lk} &\ra^{kl}_{ij}&=\rh^{ik}_{jl} &
                                \rg^{kl}_{ij}&=q^{2j-2k}\rh^{lj}_{ki}\\
  (\rc^-)^{kl}_{ij}&=(\rh^{-1})^{ji}_{lk} &(\ra^-)^{kl}_{ij}&=
     (\rh^{-1})^{ik}_{jl} & (\rg^-)^{kl}_{ij}&=q^{2j-2k}(\rh^{-1})^{lj}_{ki}.
\end{align*}  
The meaning of these variants becomes clear if one relates them to
the universal $r$-form
\begin{align}
  \rc^{kl}_{ij}&=p^{-1}\rf(S(u^i_l),S(u^j_k))
  &(\rc^-)^{kl}_{ij}&=p\overline{\rf}(S(u^j_k),S(u^i_l))\label{rcmatrizen}\\
  \ra^{kl}_{ij}&=p^{-1}\overline{\rf}(u^k_j,S(u^i_l))
  &(\ra^-)^{kl}_{ij}&=p\rf(S(u^i_l),u^k_j)\label{ramatrizen}\\
  \rg^{kl}_{ij}&=p^{-1}\overline{\rf}(S(u^j_k),u_i^l)
  &(\rg^-)^{kl}_{ij}&=p\rf(u^l_i,S(u^j_k)).\label{rgmatrizen}
\end{align}
In terms of the matrix coefficients $u^i_j$ of the
vector representation $V(\omega_1)$ the matrix coefficients
of the fundamental representation $V(\omega_s)$ are given by $s$-minors.
More explicitly let $I=\{i_1{<}\dots{<}i_s\}$ and $J=\{j_1{<}\dots{<}j_s\}$
denote subsets of $\{1,\dots,N\}$. Define
\begin{align*}
\cD^I_J:=\sum_{\sigma\in S_s}(-q)^{\ell(\sigma)}
                      u^{i_{\sigma(1)}}_{j_1}\dots u^{i_{\sigma(s)}}_{j_s}
\end{align*}                      
where $S_s$ denotes the symmetric group in $s$ elements and $\ell(\sigma)$ is
the length of the permutation $\sigma\in S_s$.
Then the matrix coefficients
\begin{align*}
  \{x_I:=\cD^{\{r+1,\dots,N\}}_I\,|\,
            I=\{i_1{<}\dots{<}i_{s}\}\subset\{1,\dots,N\}\}
\end{align*}            
form a basis of a left $\UslN$-module isomorphic to $V(\omega_s)$ with
corresponding right $\oqsl$-comodule structure
\begin{align*}
  \kopr(x_I)=x_J\otimes \cD^J_I.
\end{align*}  
As for the vector representation (\ref{R}) a solution
$(\rh^{KL}_{IJ})$ of the braid relation on
$V(\omega_s)\otimes V(\omega_s)$ is given by
\begin{align}\label{Rhmatrix}
  \rf(\cD^L_I,\cD^K_J)=p^{s^2}\rh_{IJ}^{KL}.
\end{align}  
In analogy to (\ref{rcmatrizen})-(\ref{rgmatrizen}) define
\begin{align}
  \rc^{KL}_{IJ}&=p^{-s^2}\rf(S(\cD^I_L),S(\cD^J_K))
  &(\rc^-)^{KL}_{IJ}&=p^{s^2}\overline{\rf}(S(\cD^J_K),S(\cD^I_L))
     \label{Rcmatrizen}\\
  \ra^{KL}_{IJ}&=p^{-s^2}\overline{\rf}(\cD^K_J,S(\cD^I_L))
  &(\ra^-)^{KL}_{IJ}&=p^{s^2}\rf(S(\cD^I_L),\cD^K_J)\label{Ramatrizen}\\
  \rg^{KL}_{IJ}&=p^{-s^2}\overline{\rf}(S(\cD_K^J),\cD_I^L)
  &(\rg^-)^{KL}_{IJ}&=p^{s^2}\rf(\cD^L_I,S(\cD^J_K)).\label{Rgmatrizen}
\end{align}

For further details on $q$-deformed universal
enveloping algebras $\ug$ and coordinate algebras $\oqg$ consult \cite{b-KS}.

For any subset $S$ of a complex vector space $V$ the symbol $\Lin_\C$ will
denote the complex linear span of $S$. For any linear map
$A:V^{\ot n}\rightarrow V^{\ot m}$ and some tensor power
$V^{\ot p}$, $p\ge k{+}n$ the map $\id^{\ot k}\ot A\ot \id^{\ot p-k-n}$
on $V^{\ot p}$ will be denoted by the symbol $A_{k+1,\dots,k+n}$.

\section{Quantum Grassmann manifolds}\label{grass}

\subsection{Localization and Complexification}\label{local}

Let $\mathcal X$ denote an $\ug$-module algebra without zero divisors and let
$1\in S\subset \mathcal X$ denote a left and right Ore subset, i.e.~$S$ is 
multiplicatively closed and satisfies
\begin{eqnarray*}
\forall x\in \cX, s\in S\,\quad \exists x_L,x_R\in \cX,s_L,s_R\in S\quad
 s_L x=x_Ls,\,  x s_R=s x_R.\label{ore}
\end{eqnarray*}
Then the localization ${\mathcal X}(S)$ is defined \cite{b-McCoRo}.
It is possible to
give criteria when the $\ug$-module structure of $\cX$ induces a $\ug$-module
structure on $\cX(S)$ \cite{a-LuRo}. In the example considered here the
existence of such an action of $\ug$ will be straightforward.

Let now $q\in\R$ and consider the compact real form
of some quantum group $\oqg$ with real universal $r$-form $\rf$.
Assume $\cX=\oplus_\lambda \cX_\lambda$
where $\cX_\lambda$ are irreducible type 1 representations of $\ug$,
i.e.~$\cX$ is also an ${\mathcal O}_q(G)$-comodule algebra.
Let $\overline{\cX}$ denote the complex 
conjugate vector space with the opposite multiplication $x\cdot y:=yx$. 
Then $\overline{\cX}$ is also an ${\mathcal O}_q(G)$-comodule algebra with 
comultiplication
\begin{equation*}
  \overline{\kopr}_\cX\colon x\mapsto x_{(0)}
     \otimes x_{(1)}^*\in\overline{\cX}\otimes{\mathcal O}_q(G).
\end{equation*}
Assume further that $\cX$ is given by generators and homogeneous relations.
For any homogeneous $x\in\cX$ (resp. $y\in\cbX$) let $\deg(x)$
(resp. $\deg(y)$) denote the degree, i.e.~the number of generators
occurring in each summand of $x$
(resp.~$y$). Then for each $\lambda\in \C$ the tensor product
$\cX_\C^\lambda := \cX\otimes \overline{\cX}$
obtains an algebra structure by
\begin{eqnarray}
(x\otimes y)\cdot(x'\otimes y')&=&\lambda^{\deg(y)\deg(x')}(x x_{(0)}')
\otimes (y_{(0)}y')\rf(y_{(1)},x_{(1)}')\label{lamcom}
\end{eqnarray}
and turns into an ${\mathcal O}_q(G)$-comodule algebra (\cite{b-KS}, Lemma 10.31).
The algebra $\cX_\C^\lambda$ will be called the complexification of $\cX$.
If $\lambda\in\R$ then $\cX^\lambda_\C$ is a $*$-algebra with
$(x\ot y)^*=y\ot x$.

\subsection{The homogeneous coordinate ring}
The Grassmann manifold $\Gr$ of $r$-dimensional subspaces of $\C^N$ is
the $\SL$-orbit of a highest weight vector $v\in V(\omega_r)$ in the 
projective space ${\mathbb P}(V(\omega_r))$. Its homogeneous coordinate ring
$\ovr$ is generated by the set of functions 
\begin{equation}
    \{f_\ell:=\ell(\cdot\, v)|\ell\in V(\omega_r)^*\}\label{plucker}
\end{equation}
in ${\mathcal O}(\mathrm{SL}(N))$.
This approach allows a well known analogue in the $q$-deformed setting
for arbitrary flag manifolds
\cite{b-CP94},\cite{a-LaksResh92},\cite{a-TaftTo91},\cite{a-Soib92}.
We restrict ourselves to the special case of Grassmannians.
More explicitly, if $v\in V(\omega_r)$ is a highest weight vector
then $\oqvr$ is defined to be the subalgebra generated by the
matrix coefficients (\ref{plucker}) in $\SLN$.

\begin{proposition} {\upshape (\cite{a-Soib92}, Prop.~1)}
As a $\UslN$-module algebra $\oqvr$ is isomorphic to 
the direct sum $\oplus_{k\ge0}V(k\omega_r)^*$ endowed with the Cartan 
multiplication.
\end{proposition}
Recall that the Cartan multiplication is given on homogeneous 
components by the projection
\begin{equation*}
 V(k\omega_r)^*\otimes V(m\omega_r)^*\rightarrow V((k+m)\omega_r)^*.
\end{equation*} 
The $q$-deformed coordinate algebra $\oqvr$ can further be 
described in terms of generators $V(\omega_r)^*=V(\omega_s)$, $s=N{-}r$, and
relations by 
$V(\omega_{s})\otimes V(\omega_{s})
\supset V(\lambda)=0$ if $\lambda\neq 2\omega_{s}$ \cite{a-brav94},
\cite{a-TaftTo91}. The notation $s=N{-}r$ will be used throughout this paper.

The eigenvalue of the braiding $\rh_{\omega_s,\omega_s}$ induced by the
universal $r$-form on
$V(2\omega_{s})\subset V(\omega_{s})\otimes V(\omega_{s})$ is given by
$q^{s(N-s)/N}$, compare \cite{b-KS} 8.4.3, Prop.~22.
Thus the eigenvalue of the rescaled braiding $(\rh^{KL}_{IJ})$ on
$V(2\omega_{s})\subset V(\omega_{s})\otimes V(\omega_{s})$ equals $q^{s}$.
This eigenvalue occurs on no other
irreducible subspace of $V(\omega_{s})\otimes V(\omega_{s})$.
Therefore the $q$-deformed homogeneous coordinate ring
$\oqvr$ can also be defined by generators $x_I$ and relations
\begin{align*}
  x_Ix_J=q^{-s}x_Kx_L \hat{R}^{KL}_{IJ}.
\end{align*}  
The complexification $\oqvr_\C^\lambda$ can also be given explicitly in
terms of generators and relations.

\begin{lemma}\label{relationen}
Let $(x_I)$ denote the basis of $V(\omega_{s})$ from above and let $(y_I)$
denote the dual basis of $V(\omega_{s})^*=V(\omega_{r})$.
Then the complexification $\oqvr_\C^\lambda$ is isomorphic to
\begin{align*}
  \oqvr\otimes \cO_q(\cV_s)
\end{align*}  
and a complete list of relations is given by
\begin{eqnarray}
x_I x_J&=&q^{-s}\,x_K x_L \hat{R}_{IJ}^{KL}\label{xx}\\
y_I y_J&=&q^{s}\,y_K y_L (\check{R}^-)_{IJ}^{KL}\label{yy}\\
y_I x_J&=&p^{-s^2}\lambda \,x_K y_L (\acute{R}^-)_{IJ}^{KL}.\label{xy}
\end{eqnarray}
\noindent If $\lambda=p^{s^2}q^s$ then any invariant element 
$c\in V(\omega_s)\otimes V(\omega_{r})\subset \oqvr_\C^\lambda$ is central.
\end{lemma}
\begin{proof}
To prove (\ref{xy}) note that
$\kow x_I=x_J\otimes \cD^J_I$ implies
$\kow y_I=y_J\otimes S(\cD_J^I)$ and by definition
\begin{align*}
  y_Ix_J=\lambda x_Ky_L\rf(S(\cD^I_L),\cD^K_J)=
           p^{-s^2}\lambda x_Ky_L (\ra^-)^{KL}_{IJ}.
\end{align*}           
Decomposition of the tensor product $V(\omega_s)\otimes V(\omega_r)$
yields that any invariant element $c\in V(\omega_s)\otimes V(\omega_r)$
is a complex multiple of $\sum_I x_Iy_I$.
Now the last statement follows from the relations
(\ref{Rhmatrix})--(\ref{Rgmatrizen})
and the properties of the universal $r$-form.
\end{proof}
\noindent In what follows we will only consider the case $\lambda=p^{s^2}q^s$
and drop the index $\lambda$.

\subsection{$q$-Grassmann Manifolds}\label{q-gr}

The complex Grassmann manifold $\Gr$
has the structure of a homogeneous space $\Gr\cong G/K$ with
$G=\SU$ and $K=\mathrm{S}(\mathrm{U}(r)\times \mathrm{U}(N{-}r))$.
The harmonic analysis of the square
integrable functions on $\Gr$ endowed with the Haar measure is given by
\cite{b-Helga84}, V Thm.~4.3,
\begin{equation*}
  L^2(\Gr)=\widehat{\bigoplus_{\lambda\in P^+_K} V(\lambda)}.
\end{equation*}
Here $P^+_K$ denotes the set of $K$-spherical dominant weights. In the
case of $\Gr$ the set $P^+_K$ consists of all weights of the form
\cite{b-Helga84}, V Thm.~4.1, \cite{phd-stok}, Thm.~4.2.1,
\begin{equation}
  \lambda=\sum_{i=1}^{\min(r,s)} n_i(\omega_i+\omega_{N-i}),
  \qquad n_i\in\N_0
\end{equation}  
where as before $s=N{-}r$. The direct sum 
\begin{align*}
  \mathcal{O}(\Gr)=\bigoplus_{\lambda\in P^+_K} V(\lambda)
\end{align*}  
is multiplicatively closed and will be called the coordinate algebra of
$\Gr$. As $\Gr$ is a projective variety this notion deviates from the 
classical formalism of algebraic geometry.

Consider the generators $V(\omega_{s})\subset \mathcal{O}(\plr)$ of
the homogeneous coordinate ring of $\Gr$. To any $f\in V(\omega_{s})$
associate a complex conjugate function $f^*$ defined by
\begin{align*}
  f^*(p)=\overline{f(p)}.
\end{align*}  
The vector space $\overline{V(\omega_{s})}$ of complex conjugate functions
is endowed with the complex conjugate scalar multiplication,
i.e.~$(\lambda\cdot f^*)(p)=\overline{\lambda f(p)}$.
The $\SU$-invariant scalar
product on $V(\omega_{s})$ induces an isomorphism of representations
\begin{eqnarray*}
\overline{V(\omega_{s})}&\rightarrow& V(\omega_{s})^*=V(\omega_r)\\
f^*&\mapsto&\pair{\cdot}{f}
\end{eqnarray*}
which extends to  an isomorphism of algebras
\begin{align*}
  \overline{\mathcal{O}(\plr)}\rightarrow
            \mathcal{O}({\mathcal V}_{s}).
\end{align*}            
Denote the algebra of functions on the affine cone $\plr$ generated by
$V(\omega_{s})$ and $\overline{V(\omega_{s})}$ by $\ovr_\C$.
The multiplication $\ovr\otimes\cO({\mathcal V}_s)\rightarrow \ovr_\C$
can be seen to be an isomorphism. 
Let $c\in V(\omega_s)\otimes V(\omega_{r})\subset\ovr_\C$ denote a nonzero
invariant element (which is uniquely determined up to a scalar factor)
and consider the algebra of functions on $\plr$ generated by
\begin{align*}
  V=\Lin_\C\left\{ \frac{xy}{c}\,\Big|\,x\in V(\omega_{s})\subset \ovr,
                              y\in V(\omega_{r})\subset \cO(\cV_s)
  \right\}.
\end{align*}  
The representation $V$ is isomorphic to
\begin{equation*}
  V\cong\bigoplus_{k=0}^{\min(s,r)} V(\omega_k+\omega_{N-k})
\end{equation*}  
where by abuse of notation $V(\omega_0+\omega_{N-0})$ denotes the 
trivial representation.
The elements of $V$ are homogeneous of degree $0$ and therefore induce 
functions on the Grassmann manifold. As the localization ${\ovr_\C}(c)$ has no
zero divisors the elements of $V$ generate the algebra $\ogr$. Thus 
classically $\ogr$ is isomorphic to the subalgebra of ${\ovr_\C}(c)$
generated by $V$.

This construction of $\ogr$ allows a straightforward analogue in the 
$q$-deformed setting, $q\in\R\setminus\{-1,0,1\}$. Namely consider the
complexification
\begin{equation*}
  \oqvr_\C\cong\oqvr \otimes \cO_q(\cV_s)
\end{equation*}
of $\oqvr$ and let $c\in V(\omega_s)\otimes V(\omega_{r})\subset \oqvr_\C$
denote a nonzero invariant element.
By Lemma \ref{relationen} the element $c$ is central in $\oqvr_\C$. 
The subalgebra of the localization $\oqvr_\C(c)$ generated by
\begin{equation*}
  V=\Lin_\C\left\{\frac{xy}{c}\,\Big|\,x\in V(\omega_{s})\subset \oqvr,
                             y\in V(\omega_{r})\subset \cO_q(\cV_s)
  \right\}
\end{equation*}  
will be denoted by $\oqgr$ and will be called the $q$-deformed coordinate
algebra of the Grassmann manifold. This definition of
$\oqgr$ will prove useful for our purpose
as one can use the FODC on $\oqvr$ to induce a FODC on 
$\oqvr_\C(c)$ and therefore on the $q$-deformed Grassmann manifold.
The aim of the next theorem is to make contact with the definition of
$\oqgr$ in the literature.

By construction there is an inclusion
\begin{align*}
  \oqvr&\hookrightarrow \SLN\\
  x_I&\mapsto \cD_I^{\{r{+}1,\dots,N\}}
\end{align*}  
which induces a covariant 
$\ast$-algebra homomorphism
\begin{equation*}
  i_\C:\oqvr_\C\rightarrow \SUN.
\end{equation*}  
Indeed, map $y_I$ to
$(\cD_I^{\{r{+}1,\dots,N\}})^*=S(\cD^I_{\{r{+}1,\dots,N\}})$.
It follows from (\ref{rform})
that for $A=\{r{+}1,\dots,N\}$ 
\begin{equation*}
  S(\cD^I_A)\cD_J^A=\sum_{M,N,K,L}\rg^{AA}_{MN}\cD_K^M
                    S(\cD_N^L)(\ra^-)^{KL}_{IJ}.
\end{equation*}                    
Now $(\rg)^{AA}_{MN}=q^s\delta_{MA}\delta_{NA}$ and therefore
$i_\C(x_I)$ and $i_\C(y_J)$ satisfy the relation (\ref{xy}).
Note that $i_\C(c)=1$. Therefore $i_\C$ induces a map denoted by the same
symbol
\begin{equation*}
  i_\C:\oqgr\rightarrow \SLN.
\end{equation*}  
Define $K\subset\UslN$ to be the subalgebra generated by
\begin{align}
\{E_i,F_i,K_j\,|\,i\neq r,j=1,\dots,N-1\}
\end{align}
and set $K^+=\{k\in K\,|\,\epsilon(k)=0\}$.
The following theorem implies that $i_\C$ is an isomorphism of $\oqgr$
onto the $q$-deformed coordinate algebra
\begin{align*}
  \cO_q(U/K)=\{b\in \SLN \,|\,\pair{k}{b_{(1)}}b_{(2)}=0 \quad
  \text{for all $k\in \K ^+$}\}
\end{align*}  
of the Grassmann manifold defined for instance in \cite{a-DS99}(5.1).
Here $\langle \cdot,\cdot\rangle$ denotes the pairing
$\UslN\ot \SLN\rightarrow\C$ given in \cite{b-KS}, 9.4.

\begin{theorem}\label{harman}
$\oqgr\cong\bigoplus_{\lambda\in P^+_K} V(\lambda)$.
\end{theorem}
\begin{proof}
Note first that for any $k\in\N$
\begin{eqnarray}
V(k\omega_s)\otimes V(k\omega_{r})&\cong& V((k{-}1)\omega_s)\otimes 
V((k{-}1)\omega_{r})\oplus\nonumber\\
&&\,\bigoplus_{(n_1,n_2,\dots)}
        V\left(\sum_{i=1}^{\min(r,s)}n_i(\omega_i+\omega_{N-i})\right)
                                         \label{tensInd}
\end{eqnarray}
where the summation is over all $\min(r,s)$-tuples $(n_1,n_2,\dots)$ such
that $\sum_in_i=k$.
As a vector space $\oqgr\subset\oqvr_\C(c)$ is generated by subspaces
\begin{equation*}
  W_k:=\frac{1}{c^k}V(k\omega_s)\otimes V(k\omega_{r})\subset \oqvr_\C(c)
\end{equation*}  
consisting of all $k$-fold products of elements of $V$.
In $\oqgr$ the subspace $W_{k,k{-}1}\subset W_k$ isomorphic to
\begin{equation*}
  V((k-1)\omega_s)\otimes V((k-1)\omega_{r})
\end{equation*}  
is identified with $W_{k-1}$. Indeed, if $v_i\in\oqvr_\C$ denotes a highest
weight vector of
\begin{equation*}
  V(\omega_i+\omega_{N-i})\subset V(\omega_s)\otimes V(\omega_r)\subset
  \oqvr_\C
\end{equation*}  
the element $c^\ell v_1^{n_1}\dots v_r^{n_r}$, 
$\sum n_i+\ell=k$ is up to a scalar multiple the unique highest weight vector
in $V(k\omega_s)\otimes V(k\omega_r)\subset \oqvr_\C$ of weight
$\sum n_i(\omega_i+\omega_{N-i})$.
Dividing by $c^k$ implies the claimed identification of $W_{k,k{-}1}$
with $W_{k-1}$ in $\oqgr$. Now the claim of the theorem follows from
(\ref{tensInd}).
\end{proof}
Note the inclusion $i_\C(\oqgr)\subset\cO_q(U/K)$. In fact it is known that
$\im(i_\C)$ and $\cO_q(U/K)$ are isomorphic and that 
the decomposition of $\im(i_\C)$ into irreducible modules is
the same as the decomposition of $\oqgr$ in Theorem \ref{harman}
\cite{phd-stok}. This implies the following Corollary.
\begin{cor}\label{IsoUK}
The map $i_\C:\oqgr\rightarrow \cO_q(U/K)$ is an isomorphism of left
$\UslN$-module algebras.
\end{cor}

\subsection{Generators and relations}\label{gen+rel}
To present $\oqgr$ in terms of generators and relations
consider the standard basis 
$x_I\in\Lambda^{\sr}(\C^N)=V(\omega_{\sr})$  
with multi-indices
$I=(i_1<\dots<i_{\sr})\subset\{1,\dots,N\}$. If $\e_1,\dots,\e_N$
denotes the standard basis of $V(\omega_1)=\C^N$ then
\begin{align}
  x_I=\sum_{\sigma\in S_{\sr}}(-q)^{\ell(\sigma)}\e_{\sigma(i_1)}\otimes
  \dots\otimes\e_{\sigma(i_{\sr})}\in V(\omega_{\sr})\subset
  V(\omega_1)^{\otimes \sr}.\label{x_Iq}
\end{align}  
Let further $\e_1^*,\dots,\e_N^*$ denote the dual basis of
$V(\omega_{N-1})=(\C^N)^*$ and define
\begin{align}
   y_I=\sum_{\sigma\in S_{\sr}}(-q)^{\ell(\sigma)}\e^*_{\sigma(i_{\sr})}\otimes
   \dots\otimes\e^*_{\sigma(i_{1})}\in V(\omega_{r})\subset
   V(\omega_{N-1})^{\otimes \sr}.\label{y_Iq}
\end{align}
Evaluation
\begin{align*}
  \e_{i_1}^*\ot\dots\ot\e_{i_k}^*(\e_{j_k}\ot\dots\ot\e_{j_1})=\delta_{i_1j_1}
  \dots\delta_{i_kj_k}
\end{align*}  
leads to
\begin{align*}
  y_I(x_J)=\left(\sum_{\sigma \in S_s}q^{2\ell(\sigma)}\right)\delta_{IJ}
                                        =q^{s(s{-}1)/2}[\sr]!\delta_{IJ}.
\end{align*}                                        
Thus up to normalization $(y_I)_I$ is the dual basis of $(x_I)_I$.
The element
\begin{align}
  c=\sum_I x_I\ot y_I\in V(\omega_{\sr})\ot V(\omega_r)\subset
                    \oqvr_\C\label{cdef}
\end{align} 
is invariant and the generators 
$z_{IJ}=x_Iy_J/c\in\oqgr$ fulfill the relations
\begin{align}
  \tr(z):=\sum_I z_{II}=&1\label{qspur}\\
  \sum_J q^{-\sr}z_{IM}z_{NK}\grave{R}^{MN}_{JJ}=&z_{IK}\label{qprodukt}\\
  z_{IR}z_{SL}\grave{R}^{RS}_{KJ}P^{IK}_{MN}=&0\label{qxproj}\\
  z_{IR}z_{SL}\grave{R}^{RS}_{KJ}\overline{P}^{JL}_{MN}=&0.\label{qyproj}
\end{align}
Here $P$ (resp.~$\overline{P}$) denotes any projector onto a subspace
$V(\mu)\subset V(\omega_{\sr})\otimes V(\omega_{\sr})$ such that 
$\mu\neq 2\omega_{\sr}$
(resp.~$V(\mu)\subset V(\omega_{r})\otimes V(\omega_{r})$ such that 
$\mu\neq 2\omega_{r}$).
By construction of $\oqgr$ as a subalgebra of $\oqvr_\C(c)$
it is clear that the above list of relations is complete.

To gain geometric understanding of $\oqgr$ it is useful to consider
another set of generators. Let $z_{ij}$ denote the basis element of 
\begin{equation*}
  V(0)\oplus V(\omega_1+\omega_{N-1})\subset V
\end{equation*}  
defined as the image in $\oqvr_\C(c)$ of the basis vector
$\e_i\otimes \e_j^*\in V(\omega_1)\otimes V(\omega_{N-1})$ by the map
\begin{equation}
  q^{s(s{-}3)/2}\frac{[\sr]![\sr]}{c}\quad
  \parbox{3cm}{\includegraphics*{kleinerelbild.001}}\label{zij}
\end{equation}
where the graphical calculus from Appendix \ref{graphcal} is applied.

The generator $z_{ij}$ can be written more explicitly
\begin{align*}
  z_{ij}&=q^{\sr(\sr-3)/2}[\sr]![\sr] c^{-1} \sum_{{i_2,\dots,i_{\sr}}
             \atop{{k_1,\dots,k_{\sr}}\atop{l_1,\dots,l_{\sr}}}}
   \e_{k_1}\dots\e_{k_{\sr}}A^{k_1,\dots,k_{\sr}}_{i,i_2,\dots,i_{\sr}}
   \e^*_{l_{s}}\dots\e^*_{l_1}
                    \check{A}^{l_{\sr},\dots ,l_{1}}_{i_{\sr},\dots,i_2,j}
\end{align*}
where $A$ and $\check{A}$ denote the projectors onto the invariant subspaces
$V(\omega_{\sr})$ and $V(\omega_r)$ respectively. Note that
Lemma \ref{antisymlem} implies for $m_1<m_2<\dots<m_\sr$ and arbitrary
$i_1,\dots,i_\sr$
\begin{align*}
  A_{i_1,i_2,\dots,i_\sr}^{m_1,m_2,\dots,m_s}=
  \check{A}_{i_s,\dots,i_2,i_1}^{m_s,\dots,m_2,m_1}=\frac{q^{-s(s-1)/2}}{[s]!}
  \sum_{\sigma\in S_s} (-q)^{\ell(\sigma)}\delta_{i_1,\sigma(m_1)}\dots
  \delta_{i_s,\sigma(m_s)}.
\end{align*}  
Therefore
\begin{equation*}
  z_{ij}=\frac{q^{-s(s+1)/2}}{[s-1]!c}\sum_{i_2,\dots,i_s}
       (-q)^{\ell(i,i_2,\dots,i_s)+\ell(j,i_2\dots,i_s)}
          x_{\{i,i_2,\dots,i_s\}} y_{\{j,i_2,\dots,i_s\}}
\end{equation*}          
where $x_{\{i,i_2,\dots,i_s\}}$ and $y_{\{j,i_2,\dots,i_s\}}$ are zero if two
indices coincide and $\ell(i,i_2,\dots,i_s)$ denotes the length of the
permutation
which transforms $i,i_2,\dots,i_s$ into an increasing sequence.
As noted in the last section there is a $*$-algebra homomorphism
$i_\C: \oqvr_\C\rightarrow \SUN$ which satisfies
\begin{align*}
  i_\C\left((-q)^{\ell(i_1,i_2,\dots,i_s)}x_{\{i_1,i_2,\dots,i_s\}}\right)&=
  \sum_{\sigma\in S_s}(-q)^{\ell(\sigma)}
  u^{\sigma(r+1)}_{i_1}\dots u^{\sigma(N)}_{i_s}\\
  i_\C\left((-q)^{\ell(j_1,j_2,\dots,j_s)}y_{\{j_1,j_2,\dots,j_s\}}\right)&=
  \sum_{\sigma\in S_s}(-q)^{\ell(\sigma)}
  S(u_{\sigma(N)}^{j_s})\dots S(u_{\sigma(r+1)}^{j_1})
\end{align*}
where $i_1,\dots,i_s$ and  $j_1,\dots,j_s$ are not necessarily ordered.
This implies
\begin{align}
  i_\C(z_{ij})=\sum_{k=r+1}^N q^{2k-2N-1}u^k_iS(u^j_k).\label{zuSu}
\end{align}
Note that
\begin{equation}
  \parbox{2cm}{\includegraphics*{kleinerelbild.013}}=
                  \frac{c}{q^{s(s-1)/2}[\sr]!}.\label{cqbild}
\end{equation}
The generators $z_{ij}$ will be called the little generators of $\oqgr$
while the generators $z_{IJ}$ will be called the big generators.
Given a $\UslN$-module algebra ${\mathcal O}$ and any two morphisms
$z_1,z_2\in\Hom(V(\omega_1)\otimes V(\omega_{N-1}),\mathcal{O})$
define their product 
$z_1\bullet z_2 \in\Hom(V(\omega_1)\otimes V(\omega_{N-1}),\mathcal{O})$ by
\begin{equation*}
    (z_1\bullet z_2)_{ik}=
                  \sum_n (z_1)_{ia}(z_2)_{bk}\grave{R}^{ab}_{nn}=
               q^{2N+1}\sum_j q^{-2j}(z_1)_{ij}(z_2)_{jk}.
\end{equation*}
Here as above $(z_1)_{ik}$ resp. $(z_1\bullet z_2)_{ik}$ denotes the image
of $\e_i\ot\e_k^*$. Graphically $z_1\bullet z_2$ is represented by
\begin{equation*}
\parbox{4cm}{\includegraphics*{kleinerelbild.004}}
\end{equation*}
Instead of $(z_1\bullet z_2)_{ik}$ we will also write
$\sum_j (z_1)_{ij}\bullet (z_2)_{jk}$ or just $(z_1)_{ij}\bullet (z_2)_{jk}$
which indicates that $\bullet$ is
a $q$-deformed matrix multiplication. 
\begin{proposition}\label{qkleinrel}
In $\oqgr$ the following relations hold
\begin{enumerate}
  \item[1.]
    $\sum_i z_{ii}=q^{-s}[s],\quad\sum_j z_{ij}\bullet z_{jk}=z_{ik}$,
  \item[2.]
    $q^{2b-2j}z_{ij}z_{ka}\hat{R}^{kb}_{jl}\hat{R}^{il}_{cd}=q^{2l-2i}
    z_{ci}z_{jk}\hat{R}^{jl}_{id}\hat{R}^{ab}_{kl},
    \qquad\mathrm{(Reflection \,  Equation)}$
  \item[3.]
    $z_{IJ}=(-1)^{s(s-1)/2}q^{s^2}([s]!)^{-1}\, 
            \parbox{5cm}{\includegraphics*{kleinerelbild.005}}.$  
\end{enumerate}
\end{proposition}
\begin{proof}
Recall the following relations
\begin{align}\label{ruu}
  \rh\ub\ub&=\ub\ub\rh\nonumber\\
  \rg\ub\ub^c&=\ub^c\ub\rg\\
  \rc\ub^c\ub^c&=\ub^c\ub^c\rc\nonumber 
\end{align}
which follow from the properties (\ref{rform}), (\ref{rform-}) of the
universal $r$-form. By Corollary \ref{IsoUK} it suffices to verify the above
relations for the elements $i_\C(z_{ij})\in\SLN$.
Property 1. is immediately checked.
To verify the reflection equation note first that
\begin{align*}
  J_{kl}=\begin{cases}q^{2k-2N-1}& \textrm{if $k=l>r$,}\\
                0& \mathrm{else}
         \end{cases}       
\end{align*}
is a solution of the reflection equation
\begin{align}
  J_{ij}J_{ka}\rg^{jk}_{lb}\rh^{il}_{cd}=
             J_{ci}J_{jk}\rg^{ij}_{dl}\rc^{lk}_{ba}.\label{refleq}
\end{align}             
Then the relations (\ref{ruu}) imply that $z_{ij}=u^m_iS(u_n^j)J_{mn}$ is
also a solution of (\ref{refleq}) concluding the proof of 2.

The last property will be verified by explicit calculation using the graphical
calculus in Appendix \ref{graphcal}.
\end{proof}

By the above proposition $\oqgr$ can be considered as a $q$-deformed
version of the coordinate algebra of the affine variety of projectors
onto $s$-dimensional subspaces of $\C^N$.

\subsection{Inclusion of quantum Grassmann manifolds}
Consider the surjective map of Hopf algebras
\begin{align*}
  \SUN&\rightarrow \mathcal{O}_q(\mathrm{SU}(N-1))\\
  u^i_j&\mapsto \begin{cases}
                  \epsilon(u^i_j)& \textrm{if } i=1 \textrm{ or } j=1,\\
                  u^{i-1}_{j-1}& \textrm{else. }
                \end{cases}
\end{align*}  
As the little generators $z_{ij}\in\oqgr\subset\SUN$, $i,j\ge 2$  map onto
the little generators
$z_{i-1,j-1}\in\oqgrminus\subset \mathcal{O}_q(\mathrm{SU}(N{-}1))$ one obtains
a surjection
\begin{align*}
  i^*_r:\oqgr\rightarrow \oqgrminus.
\end{align*}  
Classically this surjection corresponds to the inclusion
\begin{equation*}
  \mathrm{Gr}(r{-}1,N{-}1)\rightarrow \mathrm{Gr}(r,N).
\end{equation*}  
\begin{proposition}
   Let $\cI\subset \oqgr$, $r\ge 1$ denote the ideal generated by the set
     \begin{equation*}
       \{z_{1k},z_{l1}|\, k,l=1,\dots,N\}.
     \end{equation*}  
   Then $\ker(i^*_r)=\cI$.
\end{proposition}
\begin{proof}
  Let $\cL\subset\oqgr$ denote the ideal generated by
    \begin{equation*}
      \{z_{IJ}|\, 1\in I \textrm{ or } 1\in J\}.
    \end{equation*}  
  It follows from Proposition \ref{qkleinrel}.3 and (\ref{zuSu}) that
  $\cL\subset\cI\subset\ker(i^*_r).$
  Therefore it suffices to show that the induced surjection
    \begin{equation*}
      \oqgr/\cL\rightarrow \oqgrminus
    \end{equation*}  
  is also injective.

  For any multi-index $K=(k_1<\dots<k_s)$, $k_i\in\{2,\dots,N\}$ let 
  $K'$ denote the multi-index $(k_1{-}1,\dots,k_s{-}1)$. Note that
  \begin{align}
    \grave{R}^{KL}_{IJ}=\grave{R}^{K'L'}_{I'J'}\label{grR}\\
    \hat{R}^{KL}_{IJ}=\hat{R}^{K'L'}_{I'J'}\label{acR}
  \end{align}  
  if non of the occurring multi-indices $I,J,K,L$ contains 1.
  
  Recall that (\ref{qspur})-(\ref{qyproj}) form a complete set of
  relations for $\oqgr$. Let $\tilde{z}_{IJ}\in\oqgr$ denote the canonical
  preimage of $z_{IJ}\in\oqgrminus$ obtained by raising all
  components of the multi-indices $I$ and $J$ by one. It remains to verify
  that modulo $\cL$ the elements $\tilde{z}_{IJ}\in\oqgr$ satisfy the
  defining relations (\ref{qspur})-(\ref{qyproj}) of $\oqgrminus$ .
  This clearly holds for (\ref{qspur}).  
  As to (\ref{qprodukt}) note that $(1\in I \textrm{ or } 1\in J)$ and
  $\grave{R}^{KL}_{IJ}\neq 0$ imply $(1\in K \textrm{ or } 1\in L)$.
  Thus modulo $\cL$ summation has only to be taken over all
  $J\not\ni 1$. Now (\ref{grR}) implies that the preimages
  $\tilde{z}_{IJ}$ satisfy the defining relation (\ref{qprodukt})
  of $\oqgrminus$ modulo $\cL$.
  The desired property concerning (\ref{qxproj}) follows from
  (\ref{grR}), (\ref{acR}) and the fact that $P$ can be replaced by
  $\hat{R}{-}q^s$, similarly for (\ref{qyproj}).
\end{proof}

\section{Construction of Differential Calculus}\label{construction}

\subsection{Covariant Differential Calculus}
Let $\cX$ denote a $\C$-algebra. A first order differential calculus (FODC)
over $\cX$ is a $\cX$-bimodule $\Gamma$ together with a $\C$-linear map
\begin{equation*}
  \dif:\cX\rightarrow\Gamma
\end{equation*}
such that $\Gamma=\Lin_\C\{a\,\dif b\,c\,|\,a,b,c\in\cX\}$ and $\dif$
satisfies the Leibniz rule
\begin{align*}
  \dif(ab)&=a\,\dif b + \dif a\,b.
\end{align*}    
Let in addition $\cA$ denote a Hopf algebra and
$\Delta_\cX:\cX\rightarrow \cX\otimes\cA$ a right $\cA$-comodule algebra
structure on $\cX$.
If $\Gamma$ possesses the structure of a right $\cA$-comodule
\begin{equation*}
  \Delta_\Gamma:\Gamma\rightarrow\Gamma\ot \cA
\end{equation*}
such that
\begin{equation*}
\Delta_\Gamma(a\dif b\,c)=(\Delta_\cX a)((\dif\otimes\id)\Delta_\cX b)
                          (\Delta_\cX c)
\end{equation*}
then $\Gamma$ is called covariant. A FODC $\Gamma$ over a $*$-algebra
$\cX$ is called a $*$-calculus if there exists an involution
$*_\Gamma:\Gamma\rightarrow \Gamma$ such that
$*_\Gamma(a\dif b\,c)=c^*\dif(b^*)a^*$.
For further details of first order differential calculi consult \cite{b-KS}.

\subsection{Localization and Complexification}
Let $\Gamma$ denote an $\cX$-bimodule and let $1\in S\subset \cX$ denote
an Ore-subset as in section \ref{local}. Then
\begin{align*}
  \Gamma(S)=\cX(S)\otimes_\cX\Gamma\otimes_\cX \cX(S)
\end{align*}  
is called the localization of $\Gamma$ with respect to $S$.
By the Leibniz rule the localization of a covariant FODC 
$\mbox{d}\colon\cX\rightarrow \Gamma$ allows a uniquely
determined differential $\mbox{d}_S\colon\cX(S)\rightarrow\Gamma(S)$ such that
$\mbox{d}_S|_\cX=i_S\circ\,\mbox{d}$ where $i_S$ denotes the canonical map
$i_S\colon\Gamma\rightarrow \Gamma(S)$.

To define the complexification of a FODC assume as in section \ref{local}
that $q\in\R$ and consider the compact real form of some quantum group $\oqg$.
The complex conjugate vector space $\overline{\Gamma}$ of $\Gamma$
endowed with the opposite $\overline{\cX}$-bimodule structure
$a\cdot\gamma:=\gamma a$, $\gamma\cdot b:= b\gamma$, 
$\gamma\in\overline{\Gamma},\,a,b\in\overline{\cX}$ obtains the structure of
a covariant FODC by $\overline{\dif}=\dif$ and 
$\overline\Delta_\Gamma\colon \gamma\mapsto\gamma_{(0)}\otimes
\gamma_{(1)}^*$.
Assume again $\cX$ to be generated by elements of an irreducible type~1
representation $V$ and homogeneous relations.
Consider the $\cX^\lambda_\C$-bimodules
\begin{eqnarray*}
\Gamma^{\cX} &:=&(\cX\otimes\cbX)\otimes_\cX\Gamma\otimes_\cX(\cX
                               \otimes\cbX) \\
\Gamma^{\cbX}&:=&(\cX\otimes\cbX)\otimes_{\scriptstyle \cbX} \overline{\Gamma}
                  \otimes_{\scriptstyle \cbX}(\cX \otimes\cbX)
\end{eqnarray*}
and define $\Gamma^\lambda:=\Gamma^{\cX}\bigoplus\Gamma^{\cbX}\big/ \sim$ , 
where $\sim$ denotes the equivalence relation
\begin{align}
  y\ot\dif x+\dif y\ot x \sim \lambda^{\deg(y)\deg(x)}
  \big(x_{(0)}\ot\dif y_{(0)}+\dif x_{(0)}\ot
  y_{(0)}\big)\rf(y_{(1)},x_{(1)})\label{ydx+dyx}  
\end{align}
obtained by differentiation of (\ref{lamcom}) using the Leibniz rule.
Assume that $\Gamma$ is a graded $\cX$-bimodule and
$\deg(\dif x)=\deg(x)=1$ for $x\in V$.
For $x\in V\subset\cX,\,y\in\overline{V}\subset\cbX$ impose the additional
relation
\begin{align}\label{ydx}
  y\otimes\dif x=\lambda\dif x_{(0)}\otimes y_{(0)}\rf(y_{(1)},x_{(1)}).
\end{align}  
\begin{lemma}
  Let $\Gamma_\C^\lambda$ denote the quotient of $\Gamma^\lambda$ by
  \textup{(\ref{ydx})}, then
  \begin{equation*}
    \Gamma_\C^\lambda\cong(\overline{\cX}\otimes\Gamma)
      \oplus(\cX\otimes\overline{\Gamma}).
  \end{equation*}    
\end{lemma}
\begin{proof}
By (\ref{ydx+dyx}) the vector space $\Gamma^\lambda_\C$ is equal to the
quotient of $\Gamma^\cX\oplus\Gamma^{\cbX}$ by (\ref{ydx}) and
\begin{align}\label{dyx}
 \dif y\otimes x=\lambda  x_{(0)}\otimes \dif y_{(0)}\rf(y_{(1)},x_{(1)}).
\end{align}
Thus the claim follows from
\begin{align*}
  \Gamma^\cX/(\ref{ydx})\cong \cbX\otimes \Gamma,&&
  \Gamma^{\cbX}/(\ref{dyx})\cong \cX\otimes \overline{\Gamma}.
\end{align*}
\end{proof}  
The covariant FODC $\Gamma_\C^\lambda$ will be 
called the complexification of $\Gamma$. If $\lambda\in\R$ then
$\Gamma^\lambda_\C$ is a $*$-calculus and $(\Gamma^\cX)^*=\Gamma^{\cbX}$.

\subsection{Differential calculus over $\oqvr$}\label{classcal}

As a guiding principle in the construction of FODC over $q$-spaces one demands 
that dimensions should coincide with dimensions in the classical situation.
Recall that the classical K\"ahler differential
1-forms over the homogeneous coordinate ring $\ovr$ are given by
\begin{align*}
  \Omega^1=\big(\ovr\otimes V(\omega_s)\big)\Big/R
\end{align*}  
where $R$ denotes the $\ovr$-submodule generated by the irreducible components 
$V(\lambda)\subset \Sym^2(V(\omega_s)),\lambda\neq 2\omega_s $. As the defining
relations of $\ovr$ and $\Omega^1$ are homogeneous $\Omega^1$ is a graded 
$\ovr$-module. An element is homogeneous of degree $k$ if it can be written
as $\sum_j p_j\otimes\dif x^j$ where $p_j$ are homogeneous polynomials of
degree $k{-}1$ in the generators $x^j$ of $\ovr$.
Denote the elements of degree $k$ by $\Omega_k^1$.
\begin{lemma}\label{zerlegung}
  The homogeneous components of the K\"ahler differential
  1-forms over $\ovr$ are given by
 \begin{align}
  \Omega^1_1&\cong V(\omega_{s})\label{om1}\\
  \Omega^1_2&\cong V(2\omega_{s})\oplus\bigoplus_{k=0}^t
    V(\omega_{s-(2k+1)}+\omega_{s+(2k+1)})\label{om2}\\
  \Omega^1_{k>2}&\cong V(k\omega_{s})\oplus 
                V(\omega_{s-1}+(k-2)\omega_{s}+\omega_{s+1})\label{omk}
  \end{align}
  where $t=\min(\left[\frac{s{-}1}{2}\right],
              \left[\frac{r{-}1}{2}\right])$ in \textup{(\ref{om2})}. 
\end{lemma}

\begin{proof}
By construction the $\UcslN$-module structure of $\Omega^1_1$ and $\Omega^1_2$
is given by (\ref{om1}) and (\ref{om2}).
Indeed
\begin{align*}
  \Sym^2(V(\omega_s))=\bigoplus_{k=0}^{\min([\frac{s}{2}],[\frac{r}{2}])}
  V(\omega_{s-2k}+\omega_{s+2k}).
\end{align*}
The irreducible components isomorphic to
$V(k\omega_{s})$ and $V(\omega_{s-1}+(k-2)\omega_{s}+\omega_{s+1})$ in 
$V((k{-}1)\omega_{s})\otimes V(\omega_{s})\subset\oplr\otimes V(\omega_{s})$
are nonzero in $\Omega^1_{k>2}$ because they don't occur in 
$V((k{-}2)\omega_{s})\otimes\Sym^2(V(\omega_{s}))$.

Thus it remains to check, that the other components in
$V((k{-}1)\omega_{s})\otimes V(\omega_{s})$ can be written in terms of the
relations. This is achieved by direct computation.
\end{proof}

The aim of this subsection is to construct a graded covariant first order 
differential calculus $\Gamma$ over 
$\oqvr$ with the following properties.
\begin{enumerate}
\item[1.] As a left $\oqvr$-module $\Gamma$ is generated by the
  differentials of the generators of $\oqvr$. 
\item[2.] The FODC $\Gamma$ has the same 
$\UslN$-module structure as its classical counterpart $\Omega^1$.
\end{enumerate}
In the case $r=N/2$ to obtain uniqueness we demand instead the stronger
conditions
\begin{enumerate}
\item[1a.] $\Lin_\C\{u\,\dif v|\,u,v\in V(\omega_s)\subset\oqvr\}=
       \Lin_\C\{\dif v\, u|\,u,v\in V(\omega_s)\subset\oqvr\}.$
\item[2a.] Every homogeneous component of $\Gamma$ has the same 
$\UslN$-module structure as its classical counterpart $\Omega^1$.
\end{enumerate}
Conditions 1., 2.~and 2a.~imply that one has
to define $\Gamma_k\subset\oqvr\otimes V(\omega_s)$ by the right hand side
of the expression for $\Omega_k^1$ in Lemma \ref{zerlegung}.
It remains to show that $\Gamma=\oplus_k\Gamma_k$ can be endowed with a right 
$\oqvr$-module structure and a differential $\dif\colon\oqvr\rightarrow\Gamma$ 
such that it obtains the structure of a covariant FODC and without factoring 
by a nontrivial submodule. 

To define the right module structure of $\Gamma$ note that the Leibniz rule,
the covariance, condition 1a.~and the defining relations of $\oqvr$ imply
that for any $\lambda\neq 2\omega_s$ and any
$\sum x_i\otimes y_i\in V(\lambda)\subset V(\omega_s)\otimes V(\omega_s)$
\begin{align}\label{dxy=-xdy}
  \sum_i\dif x_i y_i=-\sum_i x_i\dif y_i.
\end{align}  
Thus the right module structure of $\Gamma$ is uniquely determined by a complex
parameter $c_\Gamma$ such that
\begin{align}
  \dif t_1 t_1=c_\Gamma t_1\dif t_1 \label{dtt}
\end{align}
for any
highest weight vector $t_1$ of $V(\omega_s)$. 
From the structure of $\UslN$ one obtains that
$t_1\otimes t_2-qt_2\otimes t_1$ is a highest weight vector of 
$V(\omega_{s-1}+\omega_{s+1})$ if $t_2=E_s(t_1)$. 
On the other hand applying $E_{s}$ twice to (\ref{dtt}) yields
\begin{eqnarray*}
q^{-1}\dif t_2 t_1+\dif t_1 t_2 &=&
             c_\Gamma(q^{-1} t_2 \dif t_1+ t_1 \dif t_2)\\
\dif t_2 t_2&=& c_\Gamma t_2 \dif t_2.
\end{eqnarray*}
Rewrite these equations as
\begin{eqnarray}
\dif t_j t_j&=&c_\Gamma t_j\dif t_j\mbox{ for } j=1,2\label{com1} \\
\dif t_2 t_1&=&\frac{c_\Gamma+1}{q^{-1}+q}t_1\dif t_2 
                 +\frac{c_\Gamma q^{-1}-q}{q^{-1}+q}t_2 \dif t_1\\
\dif t_1 t_2&=&\frac{c_\Gamma q-q^{-1}}{q^{-1}+q}t_1\dif t_2 \label{com2}
            +\frac{c_\Gamma+1}{q^{-1}+q}t_2 \dif t_1.
\end{eqnarray}
Using (\ref{com1})--(\ref{com2}) to calculate $\dif t_1(t_1t_2-qt_2t_1)$
one obtains
\begin{align*}
  0=(c_\Gamma q-q^{-1})(c_\Gamma q^{-1}-q)[t_1^2\dif t_2
                     -q t_1t_2\dif t_1].
\end{align*}                     
As by assumption the expression in the square brackets does not vanish, this
implies $c_\Gamma=q^{\pm 2}$.
Thus it is proved that there exist at most two covariant FODC satisfying
the conditions 1. and 2. above. 
It follows from the covariance and the choice of $c$ that the right 
module structure of $\Gamma$ given by (\ref{dxy=-xdy}) and (\ref{dtt})
is indeed well defined. 
\noindent
The main result is summarized in the following theorem.
\begin{theorem}\label{unicalc}
Let $r\in\{1,\dots ,N{-}1\}$. Then there exist exactly two covariant 
FODC over $\oqvr$ satisfying conditions \textup{1.}~and
\textup{2.}~\textup{(}resp.~\textup{1a.}~and \textup{2a.}~in the case
$r=N/2$\textup{)} above. 
\end{theorem}

The above theorem generalizes the classification result of \cite{a-PuWo89} for
covariant FODC on quantum vector spaces to all homogeneous coordinate rings
of quantum Grassmann manifolds.
Pusz' and Woronowicz' differential calculus on the $q$-deformed vector space
is obtained by application of the above construction 
to $O_q({\mathcal V}_{N-1})$. In the case $N=2$, $r=1$ there exist
two families of FODC labelled by a complex parameter satisfying conditions
1.~and 2.~\cite{a-PuWo89}. Thus to obtain uniqueness in this case additional
requirements like 1a.~and 2a.~are indeed necessary. 

An interesting feature of the above differential calculi over
$\oqvr$ is that the commutation relations between generators and their
differentials are in general no longer given by multiplication by the
universal $R$-matrix. Indeed, the commutation relations are given by a
covariant map $A:V(\omega_s)^{\otimes 2}\rightarrow V(\omega_s)^{\otimes 2}$
with eigenvalue $-1$ on $V(\omega_s)^{\otimes 2}\setminus\Sym^2 V(\omega_s)$. 

This apparent deficit can be dealt with considering the $\oqvr$-torsion in 
$\Gamma$.
\begin{lemma}
Let $\Gamma_{tor}$ denote the $\oqvr$-torsion submodule of $\Gamma$.
Then $\Gamma_{tor}$ is concentrated in degree $k=2$ and
\begin{equation}
\Gamma_{tor}=\bigoplus_\lambda V(\lambda)\label{tork}
\end{equation}
where summation is taken over all $\lambda=\omega_{s-m}+\omega_{s+m}$ with
$1< m \le \min(s,r)$ odd.
\end{lemma}
\begin{proof}
Indeed, it follows from (\ref{omk}) and the covariance of the multiplication
that the elements of $\Gamma_{tor}$ defined by (\ref{tork}) are
annihilated by any $x\in V(\omega_s)\subset\oqvr$.

Assume on the other hand that $x\in \oqvr$, $\gamma\in\Gam$ such that
$x\gamma=0.$ Without loss of generality assume $x$ and $\gamma$ to be
homogeneous, i.e.
\begin{align*}
  x&\in V(k\omega_s)\subset \oqvr\\
  \gamma=\gamma^1{+}\gamma^2&\in V(l\omega_s)\oplus
  V(\omega_{s-1}+(l{-}2)\omega_s+\omega_{s{+}1})\subset\Gamma.
\end{align*}
If $x$ and $\gamma$ are linear combinations of weight vectors the summands
of maximal weight with respect to the lexicographic order of the weights are
also torsion. Thus we can even assume $x$ and $\gamma$ to be weight vectors.
Choose bases $(x_i)_{i\in I_x}$ of $V(k\omega_s)$ and
$(\gamma^1_i)_{i\in I_1}$ of $V(l\omega_s)$ and
$(\gamma^2_i)_{i\in I_2}$ of $V(\omega_{s-1}+(l{-}2)\omega_s+\omega_{s+1})$
consisting of weight vectors such that $x=x_{i_x}$,
$\gamma^1=\gamma^1_{i_1}$ and $\gamma^2=\gamma^2_{i_2}$ for some
$i_x\in I_x$, $i_1\in I_1$ and $i_2\in I_2$.
Denote by $(D^x)^i_j$, $(D^1)^i_j$, $(D^2)^i_j \in \SLN$ the corresponding
matrix coefficients, i.e.
\begin{align*}
  \kopr x_j&=x_i\otimes (D^x)^i_j\\
  \kopr \gamma^1_j&=\gamma^1_i\otimes (D^1)^i_j\\
  \kopr \gamma^2_j&=\gamma^2_i\otimes (D^2)^i_j.
\end{align*}  
Now $x\gamma=0$ implies
\begin{align}\label{koprgleichung}
  0=\kopr (x\gamma)=x_i\gamma^1_j\otimes (D^x)^i_{i_x} (D^1)^j_{i_1}+
                    x_i\gamma^2_m\otimes (D^x)^i_{i_x} (D^2)^m_{i_2}.
\end{align}
As $\SLN$ is an integral domain and the weight $(k+l)\omega_s$ only appears in
the product of matrix coefficients of the first summand $\gamma^1\neq 0$
implies
$x_{\mathrm{max}}\gamma^1_{\mathrm{max}}=0$
for the highest weight vectors $x_{\mathrm{max}}\in V(k\omega_s)\subset \oqvr$
and $\gamma^1_{\mathrm{max}}\in V(l\omega_s)\subset\Gamma$. 
As $x_{\mathrm{max}}\gamma^1_{\mathrm{max}}$ is a highest weight
vector of $V((k{+}l)\omega_s)\subset \Gamma$ this is a contradiction.
Thus $\gamma^1=0$, $\gamma^2\neq 0$ and (\ref{koprgleichung}) would imply
$x_{\mathrm{max}}\gamma^2_{\mathrm{max}}=0$ where
$\gamma^2_{\mathrm{max}}\in V(\omega_{s-1}+(l{-}2)\omega_s+\omega_{s+1})\subset
\Gamma$ denotes a highest weight vector. This again leads to a contradiction.
Hence $\gamma=0$.
\end{proof}

In the construction of a covariant differential calculus on $q$-deformed 
Grassmann manifolds in subsection \ref{grcalc} a localization with respect to 
an invariant element $c\in V(\omega_s)\otimes V(\omega_s)^*$ of the
complexification of $\oqvr$ will be considered. In this localization all
$\Gamma_{tor}$ vanishes. Therefore it makes sense to divide by
$\Gamma_{tor}$ and consider 
\begin{equation*}
  \Gamma^{red}:=\Gamma/\Gamma_{tor}=\bigoplus_k \Gamma_k
\end{equation*}  
where now $\Gamma_k$ denotes the homogeneous torsion free component of
degree $k$, i.e.
\begin{align*}
  \Gamma_1\cong V(\omega_s),\hspace{1cm}
  \Gamma_{k\ge 2}\cong V(k\omega_s)\oplus
  V(\omega_{s-1}+(k-2)\omega_s+\omega_{s+1}).
\end{align*}
Note that in the case of quantum vector spaces there is no torsion and 
therefore $\Gamma^{red}=\Gamma$ is still the calculus of \cite{a-PuWo89}.

Consider the complexification $\Gamma_\C^{red}$ of the calculus
$\Gamma^{red}$ determined by the constant $c_\Gamma=q^2$.
By Lemma \ref{relationen} and the general construction a complete set of
relations of $\Gamma_\C^{red}$ is given by
\begin{align}
x_Ix_J&=q^{-s}x_K x_L \hat{R}^{KL}_{IJ}\label{xx2}\\
\dif x_Ix_J&=q^{2-s}x_K \dif x_L \hat{R}^{KL}_{IJ}\label{dxx}\\
y_Iy_J&=q^{s}y_K y_L (\check{R}^-)^{KL}_{IJ}\\
\dif y_Iy_J&=q^{s-2}y_K \dif y_L (\check{R}^{-1})^{KL}_{IJ}\label{dyy}\\
y_Ix_J&=q^sx_K y_L (\acute{R}^-)^{KL}_{IJ}\\
\dif x_Iy_J&=q^{-s} y_K \dif x_L \grave{R}^{KL}_{IJ}\\
\dif y_Ix_J&=q^s x_K \dif y_L (\acute{R}^-)^{KL}_{IJ}.
\end{align}
To obtain the relations in the case $c_\Gamma=q^{-2}$ replace (\ref{dxx}) and
(\ref{dyy}) by
\begin{align}
\dif x_Ix_J&=q^{s-2}x_K \dif x_L (\hat{R}^{-1})^{KL}_{IJ}\\
\dif y_Iy_J&=q^{2-s}y_K \dif y_L \check{R}^{KL}_{IJ}.\label{dyy2}
\end{align}  

\subsection{Differential calculus over $\oqgr$}\label{grcalc}

Classical K\"ahler differentials over a commutative algebra $A$ are given by
$\Gamma^1(A)=I/I^2$ where $I$ denotes the kernel of the multiplication
$m_A:A\otimes A \rightarrow A$. Thus if $A$ is a subalgebra of $B$ the
kernel of the induced map $\Gamma^1(A)\rightarrow \Gamma^1(B)$ is given by
\begin{align*}
  \big(\ker\, m_A\cap(\ker\, m_B)^2\big)\big/(\ker\, m_A)^2.
\end{align*}
It can be shown that in the case $A=\ogr$ and $B=\ovr_\C(c)$ this
quotient vanishes. Thus classically K\"ahler differentials over
$\ogr\subset \ovr_\C(c)$ coincide with the
$\ogr$-submodule of $\Gamma^1(\ovr_\C(c))$ generated by $\dif z_{IJ}$.
This observation and the construction of section \ref{q-gr} allow to 
introduce a $q$-deformed version of classical K\"ahler differentials 
over $\ogr$.
To do so consider the covariant FODC $\Gamma^+$ over $\cX=\oqvr$ constructed
in section \ref{classcal} uniquely determined by $c_\Gamma=q^2$.
The complex conjugate FODC $\overline{\Gamma^+}$ for 
$\overline{\cX}$ is seen to belong to $c_{\overline{\Gamma^+}}=q^{-2}$.
The complexification $\Gamma^+_\C$ of $\Gamma^+$ is a covariant FODC over
$\cX_\C$. Localize this FODC with respect to the invariant element
$c\in\cX_\C$ and denote the resulting calculus by $\Gamma^+_\C(c)$.
The canonical FODC over $\oqgr$ is defined to be the subcalculus of
$\Gamma^+_\C(c)$ generated by $\oqgr$ and is denoted by $\gamgr$.
Note that $\gamgr$ is a $*$-calculus as $\Gamma^+_\C$ is a $*$-calculus.
To perform calculations in $\gamgr$ it is
useful to point out some further relations in the localization 
$\Gamma_\C^+(c)$. Introduce partial differentiations
\begin{align}
\partial:=\dif \otimes \id:\oqvr\otimes\overline{\oqvr}\rightarrow
          \Gamma^+_\C\label{partial}\\
\overline{\partial}:=\id \otimes \dif:\oqvr\otimes\overline{\oqvr}\rightarrow 
          \Gamma^+_\C.\label{opartial}
\end{align}
Then in particular $\partial c=\dif x_Iy_I$ and 
$\overline{\partial}c=x_I\dif y_I$.
The following Lemma can be proved by direct calculation using
(\ref{xx2})--(\ref{dyy2}).

\begin{lemma}\label{Gam+rel}
In $\Gamma_\C^+(c)$ the following relations hold
\begin{align}
      &\partial c\, x_I=q^2 x_I\partial c&
                &\overline{\partial} c\, x_I=x_I\overline{\partial} c
                     \nonumber\\
      &\partial c\, y_I= y_I\partial c&
                &\overline{\partial} c\, y_I=q^{-2}y_I\overline{\partial} c
                     \nonumber\\
      &\partial c\, c^n= q^{2n}c^n\partial c&
                &\overline{\partial} c\, c^n=q^{-2n} c^n\overline{\partial} c
                     \nonumber\\
      &\partial c\, z_{IJ}= z_{IJ}\partial c&             
                &\overline{\partial} c\, z_{IJ}=z_{IJ}\overline{\partial} c
                     \nonumber\\
      &&&\nonumber\\
      & \dif x_I\, c=c \dif x_I+(q^2-1)x_I\partial c&
              &\dif y_I\, c=c \dif y_I+(q^{-2}-1)y_I\opartial c\nonumber\\
      & \dif x_I\, c^{-1}= c^{-1}\dif x_I+(q^{-2}-1)c^{-2}x_I\partial c&
    &\dif y_I\, c^{-1}= c^{-1}\dif y_I+(q^2-1)c^{-2}y_I\opartial c\nonumber
\end{align}
\begin{align}
     &\partial z_{IJ}=c^{-1} \dif x_I y_J -z_{IJ}c^{-1}\partial c&
       &\overline{\partial} z_{IJ}=c^{-1}x_I\dif y_J 
                              -z_{IJ}c^{-1}\overline{\partial} c \label{delz}\\
     &\partial z_{IJ}\,c^n=c^n\partial z_{IJ}&
       &\overline{\partial} z_{IJ}\,c^n=c^n\overline{\partial} z_{IJ}
              \nonumber\\
     &\partial z_{IJ}\,z_{KL}=q^2 z_{MN}\partial z_{OP} T^{MNOP}_{IJKL}&
       &\overline{\partial} z_{IJ}\,z_{KL}=
            q^{-2} z_{MN}\overline{\partial} z_{OP} T^{MNOP}_{IJKL}
            \label{delzz}
\end{align}     
where $T=\grave{R}_{23}\hat{R}_{12}\check{R}^-_{34}\acute{R}^-_{23}$.
\end{lemma}
\begin{cor}
As a left $\oqgr$-module $\gamgr$ is generated by the differentials 
$\dif z_{IJ}$. More explicitly
\begin{align*}
  \dif z_{IJ} z_{KL}=(q^2+q^{-2})z_{MN}\dif z_{OP}T^{MNOP}_{IJKL}- z_{MN}z_{OP}
  \dif z_{RS}T^{OPRS}_{ABDL}T^{MNAB}_{IJKC}C^{CD}
\end{align*}  
where $C^{CD}=\sum_Jq^{-s}\rg^{CD}_{JJ}$.
\end{cor}
\begin{proof}
All indices will be dropped in the calculations. Note first that by
(\ref{delzz})
\begin{eqnarray*}
z\dif z &=& (q^{-2}\partial z z+q^2\overline{\partial} zz)T^{-1}=
            (q^2\dif z z- (q^2-q^{-2})\partial z z)T^{-1}\\
zz\dif z&=& (q^{-4}\partial zzz+q^4\overline{\partial}zzz)
                                             T^{-1}_{1234}T^{-1}_{3456}=
            (q^4\dif zzz- (q^4-q^{-4})\partial zzz) T^{-1}_{1234}T^{-1}_{3456}
\end{eqnarray*}
which implies by (\ref{qprodukt})
\begin{align}
z\dif z T &= q^2\dif z z- (q^2-q^{-2})\partial z z\label{delingam}\\
zz\dif z T_{3456}T_{1234}C_{45}&=q^4\dif zz- (q^4-q^{-4})\partial zz.
\end{align}
Taking differences yields the desired relation.
\end{proof}

Note that (\ref{delingam}) implies $\partial z\in\gamgr$ and thus also
$\opartial z\in\gamgr$. Therefore $\gamgr$ can be written as a direct sum
\begin{align*}
  \gamgr=\Gamma_+\oplus \Gamma_-
\end{align*}  
where $\Gamma_+$ (resp.~$\Gamma_-$) denotes the $\oqgr$-subbimodule generated
by the set $\{\partial z_{IJ}\,|\, I,J\}$
(resp. $\{\opartial z_{IJ}\,|\, I,J\}$).
The differential $\partial$ (resp. $\opartial$)  endows $\Gamma_+$
(resp. $\Gamma_-$) with the structure of a covariant FODC. Relation
(\ref{delzz}) implies that $\Gamma_+$ (resp.~$\Gamma_-$) is even generated by
$\{\partial z_{IJ}\,|\, I,J\}$ (resp.~$\{\opartial z_{IJ}\,|\, I,J\}$) as a
left $\oqgr$-module. Writing (\ref{delzz}) graphically
\begin{align}\label{graphcom}
  \partial z\,z&=q^2\,\parbox{3cm}{\includegraphics*{graphcalc2.001}}&
  \opartial z\,z&=q^{-2}\,\parbox{3cm}{\includegraphics*{graphcalc2.002}}
\end{align}  
one sees that as a left and as a right $\oqgr$-module $\Gamma_+$
(resp.~$\Gamma_-$) is even generated by the differentials of the
little generators
$\{\partial z_{ij}\,|\,i,j\}$ (resp.~$\{\opartial z_{ij}\,|\,i,j\}$).
\begin{lemma}\label{littlebullet}
  In terms of the little generators of $\oqgr$ the following
  relations hold in $\gamgr$
\begin{align}
  z\bullet \partial z&=0& \partial z\bullet z&=\partial z\label{zdelz0}\\ 
  z\bullet \opartial z&=\opartial z& \opartial z\bullet z&=0.\label{delbzz0}
\end{align} 
\end{lemma}  
\begin{proof}
To verify (\ref{zdelz0}) note that the second identity follows from the
first one, the projector property $z\bullet z=z$ and the Leibniz rule.
The proof of $z\bullet \partial z=0$ is performed graphically.
If $\mu=q^{s(s-3)/2}[s]![s]$ then by definition of the little generators
\begin{align*}
z_{ij}\bullet \partial z_{jk}&= z_{ij}\bullet\left( \mu c^{-1}\,
              \parbox{1.5cm}{\includegraphics*{graphcalc2.003}}
              -c^{-1}z_{jk}\partial c\right)\\
            &=\frac{\mu^2}{c^2}\quad
              \parbox{4cm}{\includegraphics*{graphcalc2.004}}
              -c^{-1}z_{ik}\partial c\\
            &=\frac{\mu^2}{q^s c^2}\,
              \parbox{4cm}{\includegraphics*{graphcalc2.005}}
              -c^{-1}z_{ik}\partial c\\
            &=\frac{\mu^2}{q^s c^2}\,
              \parbox{4cm}{\includegraphics*{graphcalc2.006}}
              -c^{-1}z_{ik}\partial c
\end{align*}
The middle crossing can be resolved by means of Lemma \ref{phivarphi}.
\begin{align*}
z_{ij}\bullet \partial z_{jk}&=\frac{\mu^2}{[s]c^{2}}\,
              \parbox{4cm}{\includegraphics*{graphcalc2.007}}
              -c^{-1}z_{ik}\partial c\\
            &=\frac{\mu^2 q^s}{[s]c^{2}}\,
              \parbox{4cm}{\includegraphics*{graphcalc2.008}}
              -c^{-1}z_{ik}\partial c\\
            &=0.
\end{align*}
The conjugate relations (\ref{delbzz0}) follow from (\ref{zdelz0}) by
application of $*$ in the $*$-calculus $\Gamma^+_\C(c)$.
\end{proof}  

Let $\A$ denote a Hopf algebra and $\B\subset\A$ a right comodule subalgebra.
Assume that $\A$ is faithfully flat as $\B$-module.
Up to translation from right to left it was shown in \cite{a-herm01} that
there exists a one to one correspondence
between $\A$-covariant FODC $\Gamma$ over $\B$ and certain left ideals
$\cL\subset\B^+$, where $\B^+=\ker(\epsilon|_\B)$.
The dimension of a covariant FODC $\Gamma$ over $\B$ is defined by
\begin{align}
\dim(\Gamma)=\dim_\C(\Gamma/\Gamma \B^+)=\dim_\C(\B^+/\cL).
\end{align}
It has been proven in \cite{a-MullSch} that $\SUN$ is a faithfully flat
$\oqgr$-module.
\begin{proposition}\label{dim}
    The differential calculi $\Gamma_+$ and $\Gamma_-$ are nonisomorphic.
    Their dimensions can be estimated by $0<\dim(\Gamma_\pm)\le r(N{-}r)$. 
\end{proposition}  
\begin{proof}
  Note first that by construction of $\gamgr$ the differentials
  $\partial z_{IJ}$ do not vanish and therefore $\Gamma_+\neq 0$.
  Thus there exist $i,j$ such that $\partial z_{ij}\neq 0$.
  Lemma \ref{littlebullet} then implies that $\Gamma_+$ and $\Gamma_-$ are
  not isomorphic.

  To prove the second property note that for any differential calculus
  $\Gamma$ with corresponding left ideal $\cL$ the relation
  $\dim(\Gamma)=0$ implies $\cL=\B^+$. By the general construction of 
  differential calculi in terms of left ideals in \cite{a-herm01} this is
  equivalent to $\Gamma=0$.

  The vector space $\Gamma_-/\Gamma_-\B^+$ is generated by $N^2$ elements
  $\opartial z_{ij}$, $i,j=1,\dots,N$. Thus to verify the upper bound it
  suffices to show that $\opartial z_{ij}\in \Gamma_-\B^+$ if
  $j>r$ or $i\le r$. This follows from Lemma \ref{littlebullet}
  which implies for $j>r$
  \begin{align*}
    \opartial z_{ij}
  =q^{2N-2j+1}\epsilon(z_{jj})\opartial z_{ij}-\opartial z_{ik}\bullet z_{kj}
  =-\opartial z_{ik}\bullet z_{kj}^+\in\Gamma_-\B^+.
  \end{align*}
  For $i\le r$ consider the relation $\opartial z_{ij}=
  z_{ik}\bullet \opartial z_{kj}$. Proposition \ref{qkleinrel}.3 implies
  $z_{IJ}\in \B^+$ if one of the
  multi-indices $I$ or $J$ contains an index $\le r$. Thus the explicit
  form of the occurring $R$-matrices in (\ref{graphcom}) implies that
  $z_{IJ}\opartial z_{KL}\in \Gamma_-\B^+$ if $I$ contains an index $\le r$.
  As $\opartial z_{ij}=z_{ik}\bullet \opartial z_{kj}$ for $i\le r$ can be
  written as a linear combination of such expressions
  $z_{IJ}\opartial z_{KL}$ the proof for $\Gamma_-$ is completed.

  The estimate concerning $\Gamma_+$ is obtained similarly.
  Note first that $\partial z_{ij}=\partial z_{ik}\bullet\,z_{kj}$ implies
  $\partial z_{ij}\in\Gamma_+\B^+$ if $j\le r$.
  On the other hand applying (\ref{graphcom}) to
  $z_{ik}\bullet\partial z_{kj}=0$ one obtains the relation
  \begin{equation*}
    \parbox{2cm}{\includegraphics*{graphcalc2.009}}=0
  \end{equation*}  
  for any $i,j$. Assume $i>r$ and apply $\epsilon$ to the right $z$ factor.
  It follows from the explicit form of the occurring $R$-matrices that one
  obtains only terms of the following types:
  \begin{itemize}
    \item Multiples of $\partial z_{kk}$ , $k<r$
    \item Multiples of $\partial z_{kj}$, $r<k<i$
    \item $\partial z_{ij}$ with coefficient
      \begin{equation*}
         q^{-1}\vep(z_{ii})+(q^{-1}-q)\sum_{r<i'<i}\vep(z_{i'i'})=q^{2r-2N}.
      \end{equation*}  
  \end{itemize}  
  Thus by induction on $i$ one gets $\partial z_{ij}\in\Gamma_+\B^+$
  for all $i>r$.
\end{proof}

For the construction of $\gamgr$ one could have also started with the
FODC $\Gamma^-$ over $\cX=\oqvr$ which is uniquely determined by
$c_\Gamma=q^{-2}$. Indeed, consider $\Gamma^-_\C$ and localize with respect
to the invariant element $c$. Define complex and complex conjugate
differentiations $\partial$ and $\overline{\partial}$ as in
(\ref{partial}) and (\ref{opartial}) with $\Gamma^+_\C$ replaced by
$\Gamma^-_\C$. In analogy to Lemma \ref{Gam+rel} one obtains
\begin{lemma}\label{Gam-rel}
In $\Gamma_\C^-(c)$ the following relations hold
\begin{align}
      &\partial c\, x_I= x_I\partial c+(q^{-2}-1)c\dif x_I&
        &\overline{\partial} c\, x_I=x_I\overline{\partial} c\nonumber\\
      &\partial c\, y_I= y_I\partial c&
        &\overline{\partial} c\, y_I=y_I\overline{\partial}c+(q^2-1)c\dif y_I
        \nonumber \\
      &\partial c\, c^n= q^{-2n}c^n\partial c&
        &\overline{\partial} c\, c^n=q^{2n} c^n\overline{\partial}c\nonumber\\
      &\partial c\, z_{IJ}= z_{IJ}\partial c+(q^{-2}-1)c\partial z_{IJ}&
                &\overline{\partial} c\, z_{IJ}=z_{IJ}\overline{\partial} c
                +(q^2-1)c \overline{\partial}z_{IJ}\nonumber\\
      &&&\nonumber\\
      &\dif x_I\, c=q^{-2}c \dif x_I& &\dif y_I\, c=q^2 c \dif y_I\nonumber\\
     &&&\nonumber\\
     &\partial z_{IJ}=q^2c^{-1} \dif x_I y_J -q^2 z_{IJ}c^{-1}\partial c&
           &\overline{\partial} z_{IJ}=q^{-2} c^{-1}x_I\dif y_J 
                  -q^{-2} z_{IJ}c^{-1}\overline{\partial} c\label{delz-}\\
     &\partial z_{IJ}\,c^n=q^{-2n}c^n\partial z_{IJ}&
       &\overline{\partial} z_{IJ}\,c^n=q^{2n} c^n\overline{\partial} z_{IJ}
                                   \nonumber\\
     &\partial z_{IJ}\,z_{KL}=q^2 z_{MN}\partial z_{OP} T^{MNOP}_{IJKL}&
      & \overline{\partial} z_{IJ}\,z_{KL}=
            q^{-2}  z_{MN}\overline{\partial} z_{OP} T^{MNOP}_{IJKL}\
     \label{delzz-}
\end{align}     
where $T=\grave{R}_{23}\hat{R}_{12}\check{R}^-_{34}\acute{R}^-_{23}$.
\end{lemma}
\begin{cor}
  The $\oqgr$-subbimodule of $\Gamma^+_\C(c)$ generated by
  $\{\partial z_{IJ}\,|\, I,J\}$
  \textup{(}resp.~$\{\opartial z_{IJ}\,|\, I,J\}$, $\{\dif z_{IJ}\,|\, I,J\}$
  \textup{)} is isomorphic to the subbimodule of $\Gamma^-_\C(c)$ generated by
  $\{\partial z_{IJ}\,|\, I,J\}$
  \textup{(}resp.~$\{\opartial z_{IJ}\,|\, I,J\}$, $\{\dif z_{IJ}\,|\, I,J\}$
  \textup{)}.
\end{cor}
\begin{proof}
  Note first that the $\oqvr_\C(c)$ subbimodules of $\Gamma^+_\C(c)$ and
  $\Gamma^-_\C(c)$ generated by $\{\partial x_I\,|\,I\}$ are isomorphic as left
  $\oqvr_\C(c)$-modules.
  By (\ref{delz}) and (\ref{delz-}) there exists a left $\oqvr_\C(c)$-module
  isomorphism of these submodules mapping $\partial z_{IJ}\in\Gamma^+_\C(c)$
  to $\partial z_{IJ}\in\Gamma^-_\C(c)$. Thus the left $\oqgr$-submodules
  of $\Gamma^+_\C(c)$ and $\Gamma^-_\C(c)$ generated by
  $\{\partial z_{IJ}\,|\, I,J\}$ which by (\ref{delzz}) and (\ref{delzz-})
  coincide with the $\oqgr$-subbimodules generated by
  $\{\partial z_{IJ}\,|\, I,J\}$ are isomorphic as left modules.
  But (\ref{delzz}) and (\ref{delzz-}) also imply that the right module
  structures coincide.
  The claims concerning $\opartial$ and $\dif$ are obtained analogously.
\end{proof}  
The above isomorphisms of $\oqgr$-bimodules preserve the differentials.
Therefore a construction of $\gamgr$ starting from $\Gamma^-$ leads to the
same differential calculus.

In the case $N=2$, $r=1$ of Podle\'s' quantum $2$-sphere
$\cO_{q0}(S^2)=\oqgrez$ Proposition \ref{dim} implies that
$\dim(\gamgrez)=2$. It has been shown in \cite{a-herm01} that there exists
exactly one covariant first order $*$-calculus over $\cO_{q0}(S^2)$ of
dimension $2$ and that this calculus coincides with the calculus of
\cite{a-Po92}.

\section{Chern Classes}
Let $\cO(M)$ denote the coordinate algebra of an affine algebraic variety
$M$. Algebraic vector bundles over $M$ are in one to one correspondence
to projective modules over $\cO(M)$. A projective module $V$ can be
uniquely determined by a projector $p\in \mat_k(\cO(M))$, i.e.~by a
surjective map
\begin{equation*}
  p:\cO(M)^k\twoheadrightarrow V\subset \cO(M)^k,\quad p^2=p.
\end{equation*}  
By the general theory \cite{a-karoubi87} to any such projector is associated
the curvature 2-form 
\begin{equation*}
  R=p\dif p\dif p\in\mat_k(\Omega^2(M)).
\end{equation*}  
The differential forms $c_i\in \Omega^{2i}(M)$ defined by
\begin{equation*}
  \det(1+R)=1+c_1+c_2+\cdots+c_k
\end{equation*}  
are closed. The corresponding cohomology classes are the Chern classes of
the vector bundle. The differential forms $c_1,\dots,c_k$ can be expressed
in terms of the closed forms
\begin{align}
  \ch_l(V)=\frac{1}{l!}\tr(R^l).\label{chl}
\end{align}  
Consult the first chapter of \cite{a-karoubi87} for further details.

In the example of the Grassmann manifold $\Gr$ resp.~the homotope affine
algebraic variety with coordinate ring
$\C[z_{ij}|\,i,j=1,\dots,N]\big/(z^2=z,\tr(z)=N{-}r)$
the matrix $z$ with entries $z_{ij}$ can be considered as a projector
describing the tautological bundle over $\Gr$. The Chern classes of the tautological
bundle over $\Gr$ generate the cohomology ring $H^*(\Gr)$.

The aim of this section is to construct a $q$-deformed analogue of the
differential forms $\ch_l(V(\Gr))$ associated with the module of sections
$V(\Gr)$ of the tautological bundle. Consider the universal higher order
differential calculus $\Gamma^*_q(\Gr)$ over $\oqgr$ with first order
calculus $\gamgr$. The matrix $z=(z_{ij})\in \mat_N(\oqgr)$ of little
generators can be considered as a $q$-deformed analogue of the projector
which describes the tautological bundle. Define $q$-deformed analogues of the
differential forms (\ref{chl}) by
\begin{align}
  \ch_{q,l}(\Gr)=\frac{1}{[l]!}\tr(z\bullet
           \underbrace{\dif z\bullet \dif z\bullet\dots\bullet\dif z}_{2l})
                =\frac{1}{[l]!}
                  \tr(\underbrace{R\bullet R\bullet\dots\bullet R}_{l})
                  \nonumber
\end{align}
\begin{proposition}\label{closedprop}
  The differential forms $\ch_{q,l}$ are closed in
  {\upshape $\Gamma^*_q(\Gr)$}.
\end{proposition}

The main step of the proof is achieved by the following Lemma.

\begin{lemma}\label{mainstep}
  \begin{align*}
     \tr(z\bullet \opartial z\bullet\partial z\bullet
        \dots\bullet \opartial z \bullet \partial z\bullet \opartial z)&=0\\
     \tr(\partial z\bullet \opartial z\bullet\partial z\bullet
     \dots\bullet \opartial z \bullet \partial z\bullet z)&=0.
  \end{align*}
\end{lemma}
\begin{proof}
  The proof is performed graphically. By the commutation relations
  (\ref{delzz}) resp.~(\ref{graphcom}) one has
\begin{align*}
    \tr(z\bullet \opartial z\bullet\partial z\bullet &
        \dots\bullet \opartial z \bullet \partial z\bullet \opartial z)=\\
    =&\,q^2\,\parbox{8cm}{\includegraphics*[scale=.8]{graphchern.001}}
\end{align*}
\begin{align*}
    =&\,q^2\,\parbox{8cm}{\includegraphics*[scale=.8]{graphchern.002}}\\
    =&\,q^2 \tr(\opartial z\bullet\partial z\bullet \dots\bullet
            \opartial z \bullet \partial z\bullet \opartial z\bullet z)\\
    =&\,0        
\end{align*}
where the last equation follows from (\ref{delbzz0}).
The second relation is verified analogously.
\end{proof}  

\begin{proof}[Proof of Prop. \ref{closedprop}]
Note that in $\Gamma^*_q(\Gr)$
\begin{align}
    z\bullet\dif z\bullet \dif z = \dif z\bullet \dif z
    -\dif z\bullet z\bullet \dif z =\dif z\bullet \dif z\bullet z\label{zdd}
\end{align}
and by Lemma \ref{littlebullet}
\begin{align}
  z\bullet\underbrace{\dif z\bullet \dif z\bullet\dots\bullet\dif z}_{2l+1}&=
   z\bullet \opartial z\bullet \partial z\bullet\dots\bullet\opartial z
                        \label{zddd}\\
  \underbrace{\dif z\bullet \dif z\bullet\dots\bullet\dif z}_{2l+1}\bullet z&=
   \partial z\bullet \opartial z\bullet\dots\bullet\partial z \bullet z.
                        \label{dddz}
\end{align}
Lemma \ref{mainstep} and the relations (\ref{zdd}) - (\ref{dddz}) imply
\begin{align*}
  \dif (\ch_{q,l}(\Gr))=
    &\frac{1}{[l]!}\tr(\underbrace{\dif z\bullet \dots\bullet \dif z}_{2l+1})\\
  = &\frac{1}{[l]!}\left(\tr(z\bullet
        \underbrace{\dif z\bullet \dots\bullet \dif z}_{2l+1})
             + \tr(\dif z\bullet z\bullet
        \underbrace{\dif z\bullet \dots\bullet \dif z}_{2l})\right)
\end{align*}
\begin{align*}
  = &\frac{1}{[l]!}\left(\tr(z\bullet
        \underbrace{\dif z\bullet \dots\bullet \dif z}_{2l+1})
         +\tr(\underbrace{\dif z\bullet \dots\bullet \dif z}_{2l+1}
           \bullet z)\right)\\
  = &\frac{1}{[l]!}\Big(\tr(z\bullet
         \opartial z\bullet \partial z\bullet\dots\bullet\opartial z)
         +\tr(\partial z\bullet \opartial z\bullet\dots\bullet\partial z
              \bullet z)\Big)\\
  = & 0.          
\end{align*}  
\end{proof}

\begin{cor}
The differential forms {\upshape $\ch_{q,l}(\Gr)$} are central in
{\upshape $\Gamma^*(\Gr)$}.
\end{cor}
\begin{proof}
  It follows from (\ref{graphcom}) that the differential forms
  $\ch_{q,l}(\Gr)$ commute with the big generators $z_{IJ}$. Differentiation
  implies the claim.
\end{proof}  

\begin{appendix}

\section{Graphical calculus}  \label{graphcal}

The finite dimensional type 1 representations of the Hopf algebra $\ug$ form a
ribbon category $\ug{-}Rep$ (\cite{b-KS}, 8 Prop.~19, 21,
resp.~\cite{b-Turaev} for the notion of ribbon category).
Thus there exists a functor ${\mathcal F}_{\ug{-}Rep}$
from the category of $\ug{-}Rep$ coloured ribbon graphs to $\ug_{Rep}$
defined in \cite{b-Turaev}, Thm.~2.5, which allows a graphical calculus for
morphisms in $\ug{-}Rep$. Restrict to the
case $\gfrak=\slfrak_N$. On the second tensor power of the vector 
representation the braiding of $\ug{-}Rep$ is given by multiplication by the 
$R$-Matrix $p\rh$ (\ref{rpR}). To suppress the rational exponent
$p=q^{-1/N}$ in this action a rescaled version of the graphical calculus is
considered. 

More explicitly we restrict ourselves to tensor powers
$V_1\ot V_2\ot\dots\ot V_k$ where $V_i$, $i=1,\dots k$, denotes either
the vector representation $V$ or its dual $V^*$. With respect to the
standard basis $\e_1,\dots,\e_N$ of $V$ (resp.~$\e_1^*,\dots,\e_N^*$ of
$V^*$) crossings of lines correspond to
multiplication with rescaled $R$-matrices as follows
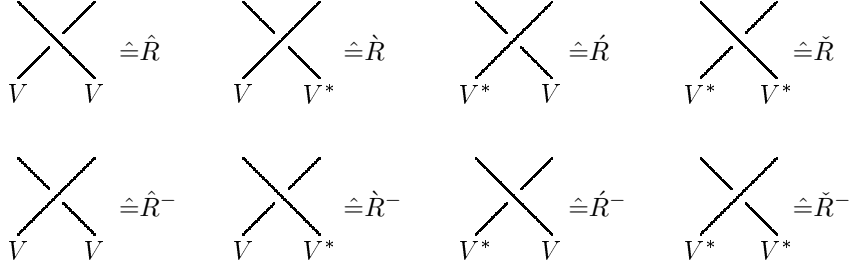
\begin{align*}
  \parbox{1.5cm}{\includegraphics*{rmat.001}}\hat{=}&\rh&
  \parbox{1.5cm}{\includegraphics*{rmat.006}}\hat{=}&\rg&
  \parbox{1.5cm}{\includegraphics*{rmat.007}}\hat{=}&\ra&
  \parbox{1.5cm}{\includegraphics*{rmat.004}}\hat{=}&\rc\\
  &&&&&&&\\
  \parbox{1.5cm}{\includegraphics*{rmat.005}}\hat{=}&\rh^-&
  \parbox{1.5cm}{\includegraphics*{rmat.002}}\hat{=}&\rg^-&
  \parbox{1.5cm}{\includegraphics*{rmat.003}}\hat{=}&\ra^-&
  \parbox{1.5cm}{\includegraphics*{rmat.008}}\hat{=}&\rc^-
\end{align*}  
and \parbox{1.5cm}{\includegraphics*{rmat.009}}
denotes the canonical inclusion $\C\rightarrow V\ot V^*$.

Projections onto sub-representations will be denoted by boxes 
containing the highest weight of the sub-representation. Thus
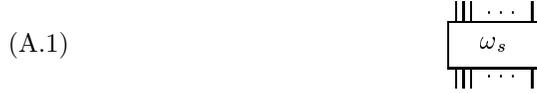
\begin{align}\label{omegarproj}
\parbox{1cm}{\includegraphics*{graphcalc3.001}}
\end{align}
denotes the projection onto the alternating sub-representations $V(\omega_s)
\subset V^{\otimes s}$. If a sub-representation occurs with some multiplicity
the box denotes the projection onto the whole isotypic component. 
Boxes labelled by $\sum_i n_i\omega_i$ will be considered as projectors of
$V^{\ot \sum in_i}$ or $(V^*)^{\ot \sum (N{-}i)n_i}$. Thus in- and outgoing
lines and lines between boxes are often dropped. Lines labelled by
$s\in \N$ represent $s$ lines. With respect to the standard basis of $V$
(resp.~$V^*$) diagrams with ingoing lines labelled by indices are identified
with the image of the corresponding basis vectors. Similarly for morphisms
on $V(\omega_s)$ (resp.~$V(\omega_s)^*$) and the standard basis $(x_I)$
(resp.~$(y_I)$) defined by (\ref{x_Iq}) (resp.~(\ref{y_Iq})).

Note that different from \cite{b-Turaev} here diagrams are considered as
morphisms from the bottom to the top. This convention is motivated by the
choice to consider right comodule algebras, i.e.~generators with lower indices.

Let $\cO$ denote a $\UslN$-module algebra with generators $(\bar{x}_I)$.
Assume that the generators $(\bar{x}_I)$ form the basis of a finite dimensional
representation $\overline{V}$ of $\UslN$. Assume in addition that there is
a given isomorphism between  $\overline{V}$ and some sub-representation of
$V^{\ot n}\ot (V^*)^{\ot m}$ for some $m,n$. Writing elements $\bar{x}$
into upper boxes with $n$ ingoing lines labelled by $V$ and $m$ ingoing lines
labelled by $V^*$ we will also consider the corresponding morphism as a map
to the algebra $\cO$. 

The following technical lemmata give a number of 
useful simplification results for morphisms in $\UslN{-}Rep$ which finally
lead to a proof of Proposition \ref{qkleinrel}.3.

\begin{lemma}\label{graph1}
\begin{equation*}
  \parbox{3cm}{\includegraphics*{graphcalc1.003}}=q^{2n_{\lambda,\mu}}
  \parbox{3cm}
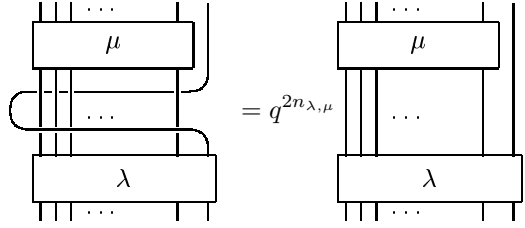
{\includegraphics*{graphcalc1.002}}
\end{equation*}  
where $n_{\lambda,\mu}=i{-}j$ if the Young tableau of $V(\lambda)$ is obtained
from the Young tableau of $V(\mu)$ by adding a box in the $(j,i)$-position,
or in terms of fundamental weights $\lambda=\mu+\omega_j-\omega_{j-1}$.
\end{lemma}

\begin{proof}
Express $\mu=\sum_k \mu_k\omega_k$, $\mu_k\in\N_0$ in terms of the fundamental
weights. By \cite{b-KS}, Prop.~8.22, the exponent $n_{\lambda,\mu}$ is given by
\begin{equation*}
  2\frac{|\mu|}{N}-(\mu,\mu+2\rho)-(\omega_1,\omega_1+2\rho)+
     (\mu+\omega_j-\omega_{j-1},\mu+\omega_j-\omega_{j-1}+2\rho),
\end{equation*}     
where $|\mu|=\sum_k k\mu_k $ denotes the number of boxes of the Young
diagram of $V(\mu)$. This expression can be simplified using the relations
\begin{align*}
(\rho,\omega_i)&=\frac{1}{2}i(N-i)\\
(\omega_i,\omega_j)&=\frac{1}{N}i(N-j)
\end{align*}
which hold for $j\ge i$.
\end{proof}

\begin{lemma}\label{phivarphi}
For $k=0,1,\dots,s$ define morphisms by
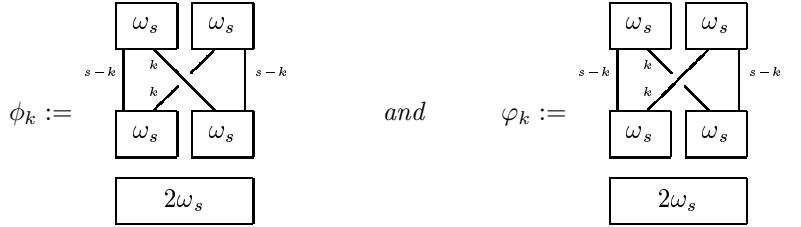
\begin{align*}
 \phi_k&:=\parbox{3cm}{\includegraphics*{graphcalc3.005}}&\textrm{and}&
 &\varphi_k&:=\parbox{3cm}{\includegraphics*{graphcalc3.004}} .
\end{align*}  
Then
\begin{align*}
     \phi_k&=\frac{q^{k(s-k+1)}}{ {s\choose k}_q }\,\phi_0 &\textsl{and}&&
  \varphi_k&=\frac{q^{-k(s-k+1)}}{ {s\choose k}_q }\,\varphi_0 .
\end{align*}  
\end{lemma}

\begin{proof}
For $k=0,1,\dots,s$ and $m=0,1,\dots,k$ define morphisms by
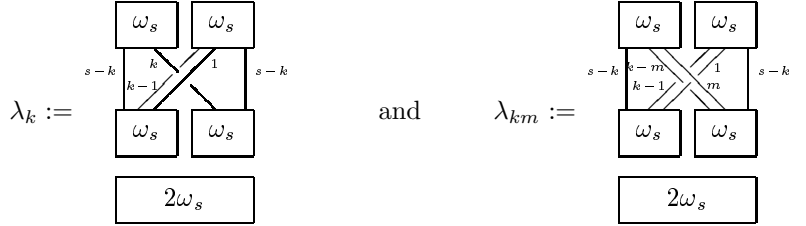
\begin{align*}
  \lambda_k:=\parbox{3cm}{\includegraphics*{graphcalc3.006}} &&\textrm{and}
  && \lambda_{km}:=\parbox{3cm}{\includegraphics*{graphcalc3.007}}.
\end{align*}
The relation $\rh^{-1}=\rh-\qd\,\id $ implies for $m<k$
\begin{align*}
  \lambda_{km}=\lambda_{k,m{+}1}-\qd q^{-2m}\phi_{k{-}1}
\end{align*}
and therefore
\begin{align}\label{lamphi}
  \lambda_k=\lambda_{k0}=\lambda_{kk}-\qd\left(\sum_{m=0}^{k-1}q^{-2m}\right)
          \phi_{k{-}1}=\phi_k-\qd q^{1-k}[k]\phi_{k-1}.
\end{align}
 On the other hand by Lemma \ref{graph1} one has
\begin{align*}
             \parbox{3cm}{\includegraphics*{graphcalc3.005}}&=
  q^{2(s{-}k)} \parbox{3cm}{\includegraphics*{graphcalc3.008}}=
  q^{2(s{-}k)} \parbox{3cm}{\includegraphics*{graphcalc3.009}}\\
  &=q^{2(s{-}k{+}1)} \parbox{3cm}{\includegraphics*{graphcalc3.006}}
\end{align*}
and therefore
\begin{align}\label{philam}
  \phi_k=q^{2(s-k+1)}\lambda_k.
\end{align}
The relations (\ref{lamphi}) and (\ref{philam}) imply
\begin{align*}
  \phi_k&=\frac{\qd q^{1-k}[k]}{1-q^{-2(s-k+1)}}\phi_{k-1}
        =\frac{q^{2(1-k)+s}[k]}{[s-k+1]}\phi_{k-1}
        =\prod_{l=1}^k\left(\frac{q^{s-2(l-1)}[l]}{[s-l+1]}\right)\phi_0\\
        &=\frac{q^{k(s-k+1)}}{ {s\choose k}_q }\,\phi_0. 
\end{align*}
The claim concerning $\varphi_k$ can be obtained from the result for
$\phi_k$ as follows.
\begin{align*}
  \varphi_k&=q^{-2k(s-k)} \parbox{3cm}{\includegraphics*{graphcalc3.010}}
            =q^{-2k(s-k)} \parbox{3cm}{\includegraphics*{graphcalc3.011}}\\
           &=q^{-2k(s-k+1)} \parbox{3cm}{\includegraphics*{graphcalc3.005}}
            =\frac{q^{-k(s-k+1)}}{ {s\choose k}_q }\,\phi_0            
\end{align*}
where in the first and in the third equation \cite{b-KS}, Prop.~8.22, is used
again.
\end{proof}

Let $A_s$ denote the antisymmetrizer, i.e.~the projector onto
$V(\omega_s)\subset V^{\otimes s}$ given graphically by (\ref{omegarproj}).
\begin{lemma}\label{antisymlem}
  \begin{align}
    A_{s+1}&=\frac{1}{[s+1]}\left(q^s(A_s)_{1,\dots,s}-[s](A_s)_{1,\dots,s}
            \rh_{s,s+1}(A_s)_{1,\dots,s}\right)\label{Ar1}\\
           &=\frac{1}{[s+1]}\left(q^s(A_s)_{2,\dots,s{+}1}
             -[s](A_s)_{2,\dots,s{+}1}\rh_{1,2}(A_s)_{2,\dots,s{+}1}\right)
             \label{Ar2},
  \end{align}
  where $(A_s)_{1,\dots,s}$ denotes the action of $A_s$ on the first $s$
  tensor components of $V^{\otimes s{+}1}$, similarly $(A_s)_{2,\dots,s{+}1}$,
  $\rh_{s,s+1}$ and $\rh_{1,2}$.
\end{lemma}

\begin{proof}
If $P_{s+1}$ is defined by the right
hand side of (\ref{Ar1}) then
\begin{align*}
  P_{s+1}|_{V(\omega_{s+1})\subset V^{\otimes (s{+}1)}}=
        \frac{1}{[s+1]}(q^s+q^{-1}[s])=1.
\end{align*}  
On the other hand it follows from Lemma \ref{phivarphi} for $k=1$ that
$P_{s+1}$ restricted to
$V(\omega_s+\omega_1)\subset V(\omega_s)\otimes V\subset V^{\otimes (s{+}1)}$
does indeed vanish.
This proves (\ref{Ar1}) and the second expression is obtained similarly.
\end{proof}

For $k=0,1,\dots,s$ and $l=0,1,\dots,s$ consider the morphisms defined by
\begin{align*}
  \Psi_{kl}&:=\parbox{3cm}{\includegraphics*[scale=.8]{graphcalc3.012}}
  &\textrm{and}&&
  \Psi_{k}&:=\parbox{3cm}{\includegraphics*[scale=.8]{graphcalc3.013}}.
\end{align*}
The irreducible representation $V(\omega_s{+}\omega_t)$ occurs in
$V(\omega_{s-k})\otimes V(\omega_{t+k})$ with multiplicity one. Thus there
exist coefficients $\psi_{kl}\in\C$ such that $\Psi_{kl}=\psi_{kl}\Psi_{k}$.
\begin{lemma}\label{psicoeff}
  The coefficients $\psi_{kl}$ satisfy the following relations:  
\begin{enumerate}
  \item $\psi_{kl}=\psi_{lk},$
  \item $\psi_{kl}=\psi_{kj}\psi_{jl}$ for all $k>j>l$,
  \item $\psi_{kk}=1,$
  \item $\psi_{k{+}1,k}=\frac{[s{-}t{-}k][k{+}1]}{[t{+}k{+}1][s{-}k]}.$ 
\end{enumerate}
\end{lemma}  
\begin{proof}
Define $W\subset V^{\otimes (s{+}t)}$ to be the isotypical component of
$V(\omega_s{+}\omega_t)$ in $V^{\otimes (s{+}t)}$. Define morphisms
$K=(A_{s{-}k}\otimes A_{t{+}k})|_W$ and $L=(A_{s{-}l}\otimes A_{t{+}l})|_W$.
Then the relation
\begin{align*}
\psi_{lk}\psi_{kl}K=\psi_{lk}KLK=KLKLK=\psi_{kl}^2K 
\end{align*}
proves the first property. Relations 2. and 3. follow immediately
from the definitions. To verify 4. define coefficients $\rho_k\in\C$,
$k=0,1,\dots,s$  by
\begin{align*}
\parbox{3cm}{\includegraphics*{graphcalc3.014}}=\rho_k\Psi_k.
\end{align*}
It follows from (\ref{Ar1}) with $s$ replaced by $s{-}k$ and from (\ref{Ar2})
with $s$ replaced by $t{+}k$ that
\begin{align*}
  \rho_k=\frac{q^{s{-}k}}{[s{-}k]}-\frac{[s{-}k{+}1]}{[s{-}k]}\psi_{k,k{-}1}
        =\frac{q^{t{+}k}}{[t{+}k]}-\frac{[t{+}k{+}1]}{[t{+}k]}\psi_{k,k{+}1}.
\end{align*}
This relation and the symmetry property $\psi_{k{+}1,k}=\psi_{k,k{+}1}$ imply
the recursion formula
\begin{align*}
  \psi_{k{+}1,k}=\frac{[t{+}k][s{-}k{+}1]}{[t{+}k{+}1][s{-}k]}\psi_{k,k{-}1}
                 +\frac{[s{-}t{-}2k]}{[t{+}k{+}1][s{-}k]}.
\end{align*}
Now property 4. follows by induction.
\end{proof}
The above Lemma can be used to calculate all coefficients $\psi_{kl}.$
In particular
\begin{align}
  \psi_{0l}=\psi_{l0}=\prod_{j=0}^{l{-}1}\psi_{j{+}1,j}=
            \prod_{j=0}^{l{-}1}\frac{[s{-}t{-}j][j{+}1]}{[t{+}j{+}1][s{-}j]}
            =\frac{{s\choose t{+}l}_q} {{s\choose l}_q {s\choose t}_q}. 
\end{align}  

\begin{lemma}\label{projinsert}
  \begin{align*}  
    \parbox{4cm}{\includegraphics*[scale=.8]{graphcalc3.015}}
   =[s]^s\prod_{j=0}^s {s\choose j}_q^{-2}
    \parbox{4cm}{\includegraphics*[scale=.8]{graphcalc3.016}}.
  \end{align*}  
\end{lemma}  
\begin{proof}
  Note first that the irreducible representation $V((s{-}1)\omega_s)$
  occurs in each of the tensor powers $V(\omega_{s{-}1})^{\otimes s}$ and
  $V(\omega_s)^{\otimes s{-}1}$ with multiplicity one.
  Therefore there
  exists a well defined coefficient $\lambda_s$ such that the
  morphism on the left hand side of the above equation is equal to
  $\lambda_s$ times the morphism on the right hand side.
  Lemma \ref{psicoeff} for $l=0$ and $k=s{-}t{-}1$ can be applied
  successively to simplify the diagram on the left hand side. This yields
  \begin{align*}
     \lambda_s=\prod_{t=0}^{s{-}1}\psi_{0,s{-}t{-}1}
              =\prod_{t=0}^{s{-}1}
                  \frac{[s]} {{s\choose t{+}1}_q {s\choose t}_q}
              =[s]^s\prod_{j=0}^{s} {s\choose j}_q^{-2}.   
  \end{align*}  
\end{proof}  

\begin{lemma}\label{straighten}
  Set $\mu_s=q^{s(s{-}1)}(-1)^{s(s-1)(s-2)/6}$ then
  \begin{align*}
    \parbox{4cm}{\includegraphics*[scale=.8]{graphcalc3.017}}
   =\mu_s\prod_{j=0}^s {s\choose j}_q^{-1}
    \parbox{4cm}{\includegraphics*[scale=.8]{graphcalc3.019}}
  \end{align*}
 where $V(\omega_s)$ is considered as a sub-representation of $V^{\otimes s}$,
 and
  \begin{align*}
    \parbox{4cm}
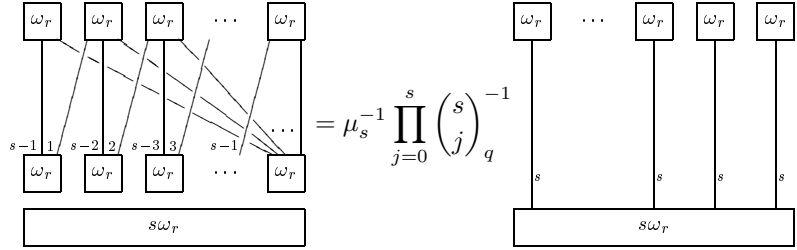
{\includegraphics*[scale=.8]{graphcalc3.018}}
   =\mu_s^{-1}\prod_{j=0}^s {s\choose j}_q^{-1}
    \parbox{4cm}{\includegraphics*[scale=.8]{graphcalc3.020}}
  \end{align*}
  where $V(\omega_r)$ is considered as a sub-representation of
  $(V^*)^{\otimes s}$.
\end{lemma}  

\begin{proof}
To prove the first relation note that Lemma \ref{phivarphi} implies
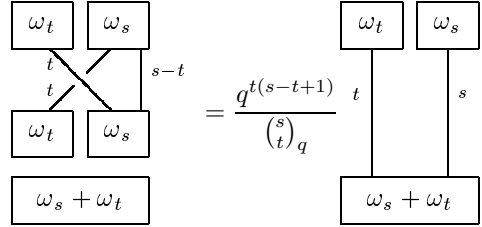
\begin{align*}
  \parbox{2.5cm}{\includegraphics*{graphcalc3.022}}
   =\frac{q^{t(s{-}t{+}1)}}{{s\choose t}_q}
    \parbox{2.5cm}{\includegraphics*{graphcalc3.023}}.
\end{align*}
Applying this successively from the right to the left hand side of the
first relation yields the desired coefficient.
The second relation is verified similarly dualizing
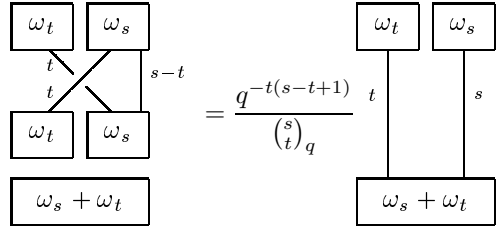
\begin{align*}
  \parbox{2.5cm}{\includegraphics*{graphcalc3.021}}
   =\frac{q^{-t(s{-}t{+}1)}}{{s\choose t}_q}
    \parbox{2.5cm}{\includegraphics*{graphcalc3.023}}.
\end{align*}
\end{proof}  
    
Now the proof of the last property of Proposition \ref{qkleinrel} can be
performed graphically. Abbreviate
$\gamma=(-q)^{-\frac{s(s{-}1)}{2}}(q^{\frac{s(s{-}3)}{2}}[s]![s])^s c^{-s}$,
then
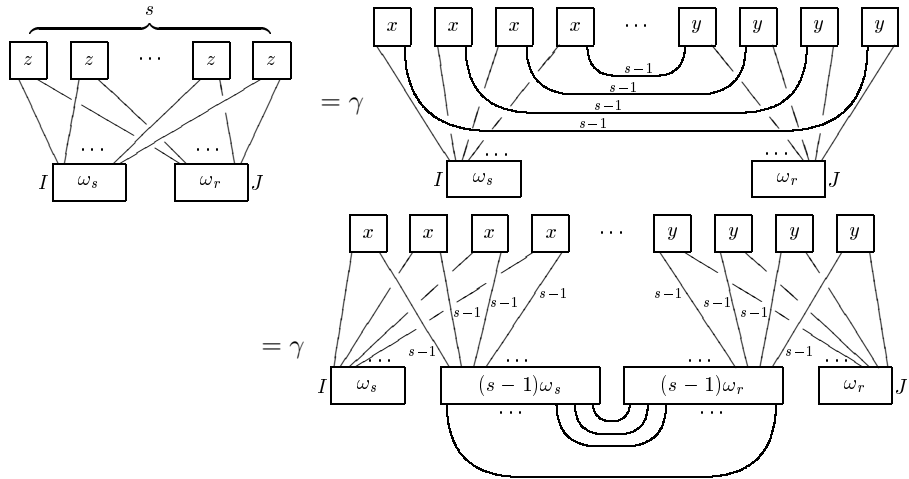
\begin{align*}
 \parbox{4cm}{\includegraphics*[scale=.8]{kleinerelbild.005}}\,&=\gamma\,\,
 \parbox{10cm}{\includegraphics*[scale=.8]{kleinerelbild.006}}\\
 =\gamma\,\,&\parbox{8cm}{\includegraphics*[scale=.8]{kleinerelbild.007}}.
\end{align*}
By Lemma \ref{projinsert} it is possible to insert $s-1$ projectors onto
$V(\omega_s)$. 
This leads to 
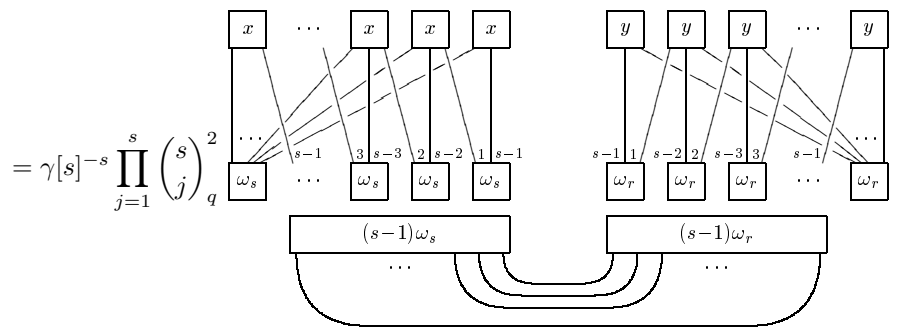
\begin{align*}
 =\gamma[s]^{-s}\prod_{j=1}^s{s \choose j}_q^2\,\,
  \parbox{9cm}{\includegraphics*[scale=.8]{kleinerelbild.012}}.
\end{align*}  
Now the left and the right half of the above diagram can be simplified using
Lemma \ref{straighten}. 
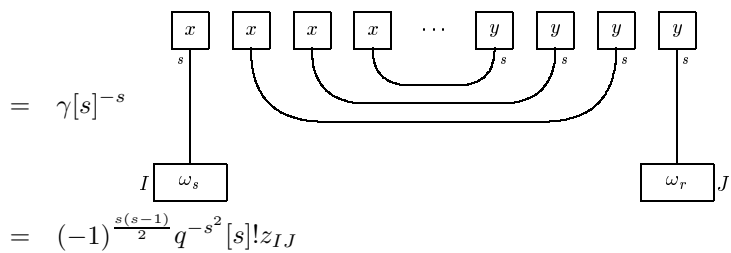
\begin{eqnarray*}
 &=&\gamma[s]^{-s}\,\,
 \parbox{9cm}{\includegraphics*[scale=.8]{kleinerelbild.011}}\\
 &=& (-1)^{\frac{s(s{-}1)}{2}}q^{-s^2}[s]! z_{IJ}
\end{eqnarray*}
which proves the last property of Proposition \ref{qkleinrel}.

\end{appendix}

\section*{Acknowledgements}
I am very indebted to Prof.~L.L.~Vaksman who suggested this work to me and
to I. Heckenberger for many hints and detailed discussions.
I am also grateful to  Prof.~K.~Schm\"udgen and U. Hermisson for helpful
comments.

\bibliographystyle{amsalpha}
\bibliography{litbank}

\end{document}